\documentclass{gtart}

\def\ifplaintex{\expandafter\ifx\csname documentclass\endcsname\relax}

\ifplaintex 
\else
\fi

\expandafter\ifx\csname beginpicture\endcsname\relax
\expandafter\ifx\csname documentclass\endcsname\relax
\input pictex \else
\input prepictex \input pictex \input postpictex \fi\fi

\def\gt{{\mathsurround=0pt\it $\cal G\mskip-2mu$eometry \&\ 
$\cal T\!\!$opology}}        

\def\gtp{{\mathsurround=0pt\it $\cal G\mskip-2mu$eometry \&\ 
$\cal T\!\!$opology $\cal P\!$ublications}}  

\def\lognumber#1{\def\thelognumber{#1}}
\def\volumenumber#1{\def\thevolumenumber{#1}}
\def\papernumber#1{\def\thepapernumber{#1}}
\def\volumeyear#1{\def\thevolumeyear{#1}}

\def\pagenumbers#1#2{\def\startpage{#1}\def\finishpage{#2}}
\def\published#1{\def\publishdate{#1}}
\def\proposed#1{\def\theproposer{#1}}
\def\seconded#1{\def\theseconders{#1}}
\def\received#1{\def\receiveddate{#1}}
\def\revised#1{\def\reviseddate{#1}}
\def\accepted#1{\def\accepteddate{#1}}

\def\asciiaddress#1{\def\theasciiaddress{#1}}

\long\def\asciiabstract#1{\long\def\theasciiabstract{#1}}

\let\\\par\let\thelognumber\relax
\let\thevolumenumber\relax\let\thepapernumber\relax
\let\thevolumeyear\relax\let\thesamplenumber\relax\let\startpage\relax
\let\finishpage\relax\let\publishdate\relax\let\receiveddate\relax
\let\reviseddate\relax\let\accepteddate\relax\let\theasciititle\relax
\let\theasciiauthors\relax\let\theasciiaddress\relax
\let\theasciiabstract\relax
\let\theasciiemail\relax\let\theshortauthors\relax\let\theshorttitle\relax

\long\def\maketitlep{   

\count0=\startpage

\gt\hfill      
\beginpicture
\setcoordinatesystem units <0.33truein, 0.33truein> point at 2.2 0.9
\setplotsymbol ({$\cal G$})
\plotsymbolspacing=9truept
\circulararc 315 degrees from 0 1 center at 0 0
\setplotsymbol ({$\cal T$})
\circulararc 315 degrees from 1 -1 center at 1 0
\endpicture
\break
{\small\ifx\thesamplenumber\relax 
Volume \else Sample
\fi\thevolumenumber\ (\thevolumeyear)
\startpage--\finishpage\nl
Published: \publishdate}
\vglue 0.5truein plus 0.4fil minus 0.1truein

{\parskip=0pt\leftskip 0pt plus 1fil\def\\{\par\smallskip}{\ifplaintex\large
\else\Large\fi\bf\thetitle}\par\medskip}   

\vglue 0pt plus 0.1fil 

{\parskip=0pt\leftskip 0pt plus 1fil\def\\{\par}{\sc\theauthors}
\par\medskip}

\vglue 0pt plus 0.1fil 

{\small\parskip=0pt\let\newline\\
{\leftskip 0pt plus 1fil\def\\{\par}{\sl\theaddress}\par}
\expandafter\ifx\theemail\relax    
\relax\else\vglue 5pt plus 0.02fil minus 2pt\def\\{\stdspace{\rm 
and}\stdspace} 
\cl{Email:\stdspace\tt\theemail}\fi
\ifx\theurl\relax                  
\relax\else\vglue 5pt plus 0.02fil minus 2pt\def\\{\stdspace{\rm 
and}\stdspace}
\cl{URL:\stdspace\tt\theurl}\fi\par}

\vglue 7pt plus 0.3fil minus 3pt

{\bf Abstract}
\vglue 5pt plus 0.1fil minus 2pt

\theabstract

\vglue 7pt plus 0.3fil minus 3pt

{\bf AMS Classification numbers}\quad Primary:\quad \theprimaryclass

Secondary:\quad \thesecondaryclass

\vglue 5pt plus 0.3fil minus 2pt

{\bf Keywords}\quad \thekeywords

\vglue 10pt plus 0.5fil minus 5pt

{\small  Proposed: \theproposer\hfill Received: \receiveddate\nl
Seconded: \theseconders\hfill 
\ifx\reviseddate\relax                         
Accepted: \accepteddate                        
\else
Revised: \reviseddate                          
\fi}
\eject
}       

\let\maketitlepage\maketitlep
\let\maketitle\maketitlepage

\font\phead=cmsl9 scaled 950
\font=cmsl9 scaled 1050
\font\pnum=cmbx10 scaled 913
\font\lnum=cmbx10 
\font\pfoot=cmsl9 scaled 950
\font=cmsl9 scaled 1050
\ifplaintex
\headline{\vbox to 0pt{\vskip -4.5mm\line{\small\phead\ifnum
\count0=\startpage ISSN 1364-0380 (on line)
1465-3060 (printed) \hfill {\pnum\folio}\else\ifodd\count0\def\\{ }%
\ifx\theshorttitle\relax\thetitle\else\theshorttitle\fi\hfill{\pnum\folio}
\else\def\\{ and }{\pnum\folio}\hfill\ifx\theshortauthors\relax\theauthors
\else\theshortauthors\fi\fi\fi}\vss}}
\footline{\vbox to 0pt{\vglue 0mm\line{\small\pfoot\ifnum\count0=\startpage
\copyright\ \gtp\hfill\else
\gt, Volume \thevolumenumber\ (\thevolumeyear)\hfill\fi}\vss
}}
\else
\makeatletter
\def\@oddhead{{\small\lhead\ifnum\count0=\startpage ISSN 1364-0380 (on line)
1465-3060 (printed) \hfill {\lnum\number\count0}\else\ifodd\count0
\def\\{ }\ifx\theshorttitle\relax \thetitle \else\theshorttitle\fi\hfill
{\lnum\number\count0}\else\def\\{ and }{\lnum\number\count0}
\hfill\ifx\theshortauthors\relax 
\theauthors\else\theshortauthors\fi\fi\fi}}\def\@evenhead{\@oddhead}
\def\@oddfoot{\small\lfoot\ifnum\count0=\startpage\copyright\ \gtp\hfill\else
\gt, Volume \thevolumenumber\ (\thevolumeyear)\hfill\fi}
\def\@evenfoot{\@oddfoot}
\makeatother
\fi

\newwrite\gtoutfile
\long\gdef\makeheadfile{  
{\def\\{, }\def\s{ }
\immediate\openout\gtoutfile head.xxx
\immediate\write\gtoutfile{Proxy-for: \ifx\theasciiauthors\relax
\theauthors\else\theasciiauthors\fi\s<\ifx\theasciiemail\relax\theemail\else\theasciiemail\fi>}
\immediate\write\gtoutfile{\noexpand\\}
\immediate\write\gtoutfile{Authors: \ifx\theasciiauthors\relax
\theauthors\else\theasciiauthors\fi}
{\def\\{ }\immediate\write\gtoutfile{Title: \ifx\theasciititle\relax
\thetitle\else\theasciititle\fi}}
\immediate\write\gtoutfile{Subj-class: GT or SG or MG etc}
\immediate\write\gtoutfile{MSC-class: \theprimaryclass\ifx\thesecondaryclass\relax\else, \thesecondaryclass\fi}
\immediate\write\gtoutfile{Journal-ref: Geom. Topol. \thevolumenumber
(\thevolumeyear) \startpage-\finishpage}
\immediate\write\gtoutfile{Comments: Published by Geometry and Topology at}
\immediate\write\gtoutfile{\s\s http://www.maths.warwick.ac.uk/gt/GTVol\thevolumenumber/paper\thepapernumber.abs.html}
\immediate\write\gtoutfile{\noexpand\\}
\immediate\write\gtoutfile{}
\ifx\theasciiabstract\relax
\immediate\write\gtoutfile{\theabstract}\else
\immediate\write\gtoutfile{\theasciiabstract}\fi
\immediate\write\gtoutfile{}
\immediate\write\gtoutfile{\noexpand\\}
\immediate\write\gtoutfile{}
\immediate\closeout\gtoutfile}}  

\def\maketitlepage{\maketitlep\makeheadfile}
\let\maketitle\maketitlepage

\lognumber{311}
\volumenumber{7}\papernumber{20}\volumeyear{2003}
\pagenumbers{713}{756}
\received{28 March 2003}
\revised{26 October 2003}
\published{30 October 2003}
\accepted{29 October 2003}
\proposed{Benson Farb}
\seconded{Leonid Polterovich, Joan Birman}

\usepackage{amsmath,amssymb}

\newtheorem{thm}{Theorem}[section]
\newtheorem{lemma}[thm]{Lemma}
\newtheorem{cor}[thm]{Corollary}

\newtheorem{proposition}[thm]{Proposition}

\newtheorem{prop}[thm]{Proposition}

\theoremstyle{definition}
\newtheorem{defn}[thm]{Definition}

\newtheorem{remark}[thm]{Remark}

\DeclareMathOperator{\Fix}{Fix}

\DeclareMathOperator{\Diff}{Diff}
\DeclareMathOperator{\Int}{Int}
\DeclareMathOperator{\Per}{Per}
\DeclareMathOperator{\area}{area}
\newcommand{\R}{{\mathbb R}}

\newcommand{\N}{{\mathbb N}}
\newcommand{\F}{{\cal F}}
\newcommand{\M}{{\cal M}}
\newcommand{\cR}{{\cal R}}
\newcommand{\X}{\mathfrak{X}}
\newcommand{\Z}{\mathbb Z}
\def\ti{\tilde}
\def\sinfty{S_{\infty}}
\def\ti{\tilde}
\def\rtwo{\mathbb R^2}
\def\rone{\mathbb R}
\def\tag{translation arc geodesic}
\def\nc{near cycle}
\def\calo{{\cal O}}
\def\calb{{\cal B}}
\def\G{{\cal G}}
\def\H{{\cal H}}
\def\A{{\mathbb A}}
\def\AO{{\mathbb A}^{\circ}}
\def\B{{\mathbb B}}
\def\Q{{\mathbb Q}}
\def\RH{RH(W,\partial_+ W)}

\begin{document}
\title{Periodic points of Hamiltonian\\surface diffeomorphisms}
\author{John Franks\\Michael Handel}
\address{Department of Mathematics, Northwestern 
University\\Evanston, IL 60208-2730, USA}

\secondaddress{Department of Mathematics, CUNY, Lehman 
College\\Bronx, NY 10468, USA}

\asciiaddress{Department of Mathematics, Northwestern 
University\\Evanston, IL 60208-2730, USA\\and\\Department 
of Mathematics, CUNY, Lehman 
College\\Bronx, NY 10468, USA}

\email{john@math.northwestern.edu, handel@g230.lehman.cuny.edu}

\begin{abstract}
The main result of this paper is that every non-trivial Hamiltonian
diffeomorphism of a closed oriented surface of genus at least one has
periodic points of arbitrarily high period.  The same result is true for $S^2$
provided the diffeomorphism has at least three fixed points.  In
addition we show that up to isotopy relative to its fixed point set,
every orientation preserving diffeomorphism $F \co  S \to S$ of a closed
orientable surface has a normal form.  If the fixed point set is
finite this is just the Thurston normal form.
\end{abstract}
\asciiabstract{%
The main result of this paper is that every non-trivial Hamiltonian
diffeomorphism of a closed oriented surface of genus at least one has
periodic points of arbitrarily high period.  The same result is true
for S^2 provided the diffeomorphism has at least three fixed points.
In addition we show that up to isotopy relative to its fixed point
set, every orientation preserving diffeomorphism F: S --> S of a
closed orientable surface has a normal form.  If the fixed point set
is finite this is just the Thurston normal form.}

\keywords{Hamiltonian diffeomorphism, periodic points, geodesic
lamination}

\primaryclass{37J10}
\secondaryclass{37E30}

\maketitlepage

\section{Introduction}

The main result of this paper is that every non-trivial Hamiltonian
diffeomorphism of a closed oriented surface of genus at least one has
periodic points of arbitrarily high period.  The same result is true for $S^2$
provided the diffeomorphism has at least three fixed points.  It
was previously known (see \cite{Fr3}, \cite{Fr4}) that any non-trivial area preserving homeomorphism of $S^2$ with at least three fixed points   has infinitely many periodic points.
The existence of unbounded periods is a substantially stronger
conclusion that is not only interesting in its own right but has
applications \cite{fh:sl3} to the algebraic properties of the group of
area preserving diffeomorphisms of closed oriented surfaces.

\begin{thm}  \label{thm: periodic point} 
Suppose $F\co  S \to S$ is a non-trivial, Hamiltonian
diffeomorphism of a closed oriented surface $S$ and that if $S = S^2$ then 
$\Fix(F)$ contains at least three
points. Then $F$ has periodic points of arbitrarily high period.
\end{thm}

Our second result is closely related to, and depends on, a celebrated
result of Thurston \cite{Th} which provides a canonical form, up to
isotopy relative to a finite invariant set, for any homeomorphism of a
compact surface.  With the added hypothesis that the homeomorphism is
a diffeomorphism, we obtain an analogous result for isotopy relative
to the fixed point set of the map, even if that fixed point set is 
infinite.  

\begin{thm}\label{canonical form} 
Every orientation preserving diffeomorphism $F \co  S \to S$ of a closed
orientable surface has a normal form up to isotopy relative to its
fixed point set.
\end{thm}

More details, including the definition of {\em normal form} is
provided below in \S \ref{normal_form}.

The main technical work in the paper (sections 8 through 11) is an
analysis of the dynamics of an element $f \in {\cal P}(S)$, the set of
diffeomorphism of the surface $S$ that are isotopic to the identity
relative to the set of periodic points of $f$.  (The standard example
of an element of ${\cal P}(S)$ is the time one map of a flow whose
periodic points are all fixed.)  Preservation of area plays no role in
this analysis.  But rather striking conclusions follow when the area
preserving hypothesis is added.  Indeed we ultimately prove, for
example, that the identity is the only area preserving element of
${\cal P}(S^2)$ that fixes at least three points.  One can construct
non-trivial area preserving examples in ${\cal P}(T^2)$ by stopping an
irrational flow at a closed set of points.  These examples have
non-trivial mean rotation vector and so are not Hamiltonian.  (See
Corollary~\ref{two definitions are equal}.) All elements of ${\cal
P}(S)$ have zero entropy \cite{katok:entropy}. While the converse is
not true, ${\cal P}(S)$ is a rich family of zero entropy
diffeomorphisms and the techniques that we develop here should be
useful in understanding the extent to which arbitrary zero entropy
diffeomorphisms behave like the time one maps of flows.  See also
\cite{han:ze}.

As a consequence of Theorem~\ref{thm: periodic point} and the 
techniques used to prove it we will establish the following result
which is of some independent interest.

\begin{thm}  \label{thm: non-isotopic} 
Suppose $F\co S \to S$ is a non-trivial, Hamiltonian
diffeomorphism of a closed oriented surface $S$ and that if $S = S^2$ then 
$\Fix(F)$ contains at least three
points. Then there exist $n>0$ and a finite set $P \subset \Fix(F^n)$
such that $F^n$ is not isotopic to the identity relative to $P$.
\end{thm}

We will show that a consequence of this is the following slight 
strengthening of Theorem~\ref{thm: periodic point}.

\begin{thm}  \label{thm: all high periods} 
Suppose $F\co  S \to S$ is a non-trivial, Hamiltonian
diffeomorphism of a closed oriented surface $S$ and that if $S = S^2$ then 
$\Fix(F)$ contains at least three
points. Then there exist $n>0$ and $p >0$ such that $F^n$ has a point
of period $k$ for every $k \ge p.$
\end{thm}

\section{The flux homomorphism and rotation vectors}

Suppose $S$ is a closed oriented surface and $\omega$ is a smooth
volume form.  We will generally assume a fixed choice of $\omega$ and
refer to a diffeomorphism $F\co  S \to S$ which preserves $\omega$ as an
{\em area preserving diffeomorphism}.  We denote the group of
diffeomorphisms preserving $\omega$ by $\Diff_\omega(S)$ and its
identity component by $\Diff_\omega(S)_0.$ Similarly, the identity component of $\Diff(S),$ is denoted $\Diff(S)_0$.

\begin{remark} \label{isotopy through area preserving} The inclusion of
$\Diff_\omega(S)_0$ into $\Diff(S)_0$ is a homotopy equivalence for all $S$ (see Corollary 1.5.4 of
\cite{Ban}).   It follows that an $\omega$ preserving diffeomorphism $F_1$ that is isotopic to the identity is isotopic to the identity through $\omega$ preserving diffeomorphisms.  Moreover, any   isotopy  $F_t$ of $F_1$ to the identity is homotopic relative to its  endpoints to an isotopy of $F_1$ to the identity through $\omega$  preserving diffeomorphisms.
\end{remark}

\begin{defn}  
Consider the set $K$ of paths
$F_t$ in $\Diff_\omega(S)_0$, with $t \in [0,1]$, such that
$F_0 = id$.
We let $\G_\omega(S)$ denote the set of equivalence classes in $K$ where
two paths are equivalent if they are homotopic relative to their
endpoints.  There is a group multiplication on $\G_\omega(S)$ defined
by composition, i.e., $[G_t][F_t] = [H_t]$ where $H_t$ is 
defined by $H_t(x) = G_t(F_t(x))$ for $x \in S.$

\end{defn}

\begin{remark} 
If $S$ has genus greater than $1$ then $\Diff(S)_0$ is simply
connected (see Theorem 1 of \cite{EE} or \cite{Ham2}) so $[F_t]$
depends only on $F_1$.  This provides a canonical identification of
$\G_\omega(S)$ with $\Diff_\omega(S)_0$ by identifying $[F_t]$ with
$F_1$.
\end{remark}

\begin{remark} Each $[F_t] \in
\G_{\omega}(S)$ determines a preferred lift $\ti F_1\co  \ti S \to \ti S$
of $F_1$ to the universal cover of $S$.  This is the lift of $F_1$
obtained by lifting the isotopy $F_t$ starting at the identity of $\ti
S.$ It is clearly independent of the choice of representative in the
homotopy class $[F_t]$.  We are primarily concerned here with the case
that the genus is $1$.  In this case (see Theorem 1 of \cite{EE} or \cite{Ham1}) the lifting provides a canonical identification of
$\G_\omega(T^2)$ with the group of diffeomorphisms of $\R^2$ which are
lifts of area preserving diffeomorphisms of $T^2$ isotopic to the
identity. In general this is the viewpoint we will most often take on
$\G_\omega(T^2)$.

\end{remark}

We want now to define the {\it flux homomorphism} from 
$\G_\omega(S)$ to $H^1(S,R).$  We will largely follow
Chapter 3 of \cite{Ban}, where it is shown that
for any smooth $\sigma\co  S^1 \to S$ and any path $G_t$ in 
$\Diff_\omega(S)_0$ starting at the identity, the map 
$G_\sigma \co  S^1 \times I \to S$ defined by 
$G_\sigma(x,t) = G_t(\sigma(x))$ has the property that
the value of 
\[
\int_{S^1 \times I} G_\sigma^*(\omega) 
\]
depends only on the equivalence class $[G_t]$ of $G_t$ in
$\G_\omega(S)$ and the homology class of the singular cycle
$[\sigma].$   Moreover, the function assigning
the value to the homology class $[\sigma]$ is a homomorphism.   
These properties are summarized in the following result:

\begin{thm}[See \cite{Ban}]\label{thm: flux homo}
There is a unique homomorphism $$\F\co  \G_\omega(S) \to H^1(S,R)$$ with
the property that 
\[
\F([G_t])([\sigma]) = \int_{S^1 \times I} G_\sigma^*(\omega) 
\]
for any $[G_t] \in \G_\omega(S)$ and any smooth singular cycle 
$\sigma \co  S^1 \to S.$ The homomorphism $\F$ is called the {\em flux homomorphism.}
\end{thm}

\begin{remark}\label{remark: flux area}
 The flux is well defined when
$G_t$ is an isotopy from the identity to $G_1 \in \Diff_{\omega}S$
even if $G_t$ is not area preserving for $t \in (0,1).$  Indeed, 
Remark~\ref{isotopy through area preserving} implies that   
 the path 
$G_t$ in $\Diff(S)_0$ is homotopic relative to its endpoints
to a path $H_t \in \Diff_{\omega}S \subset \Diff(S)_0$ and 
  Stokes' theorem implies that the integrals defining
$\F(H_t)$ and $\F(G_t)$ have the same value.
\end{remark}

\begin{remark}\label{remark: flux cover}
The value of the integral defining $\F([G_t])([\sigma])$ is the same as the
value of a similar integral in a covering space of $S$ and we will
want to make use of this.
More precisely, suppose $p\co  \bar S \to S$ is a
covering space, $\sigma\co  S^1 \to S$ is a singular cycle and that
 $\bar \sigma$ is a lift of $\sigma$ to $\bar S$. Let
$\bar \omega = p^*(\omega)$, and let
$\bar G_t$ be the lift of $G_t$ with $\bar G_0 = id.$ 
We observe that
\[
\F([G_t])([\sigma]) = \int_{S^1 \times I} \bar G_{\bar \sigma}^*(\bar \omega). 
\]
This is because $p \circ \bar G_{\bar \sigma} = G_{\sigma}$ and hence
$G_\sigma^*(\omega) = \bar G^*_{\bar \sigma}(p^*(\omega)) = 
\bar G^*_{\bar \sigma}(\bar \omega).$
\end{remark}

Closely related to the flux homomorphism is the {\em mean rotation
vector} associated to elements of $\G_{\omega}(S)$ which we now define
following \cite{Fr5}.  We want to define a ``homological
displacement'' of a point $x \in S$ by $[G_t]$.  We will do this by
considering the path $G_t(x)$ from $x$ to $G_1(x)$ and closing it to
form a loop whose homology class will define this displacement.  There are
choices of how to close it.  We consider $S$ as a convex polygon with
edges identified.  We choose a subset $X$ of this polygon consisting
of its interior and a subset of its boundary chosen so the
identification map to $S$ is a bijection.  Then each path in $S$ can
be closed by taking the points of $X$ corresponding to its endpoints
and joining them to a base point in the interior of the polygon (see
\cite{Fr5} for more details).

It is easy to see that if the subset of the boundary of the polygon
which is in $X$ is a Borel measurable set then the function $\delta\co  S
\to H_1(S,R),$ which assigns to $x$ the homology class of the loop
formed by closing the path $G_t(x)$ in this way, is a bounded Borel
measurable function.  We will call $\delta$ a {\em homological
displacement function} associated to $[G_t].$ 

Concatenating the paths $G_t(G_1^i(x))$ for $0 \le i < n$ gives a path
from $x$ to $G_1^n(x)$ and the homology class of the loop formed by
closing this path is $\sum_{i=0}^{n-1} \delta(G_1^i(x)),$ (See
\cite{Fr5} for more details).

\begin{defn}\label{def: point rotation vec} Suppose that $\delta\co  S \to H_1(S,R)$ is a homological
displacement function associated to an isotopy $[G_t].$
We define the {\em rotation vector} of a point $x \in S$ for the
isotopy $G_t$ to be
\[
\rho(x, G_t) 
= \lim_{n \to \infty} \frac{1}{n} \sum_{i=0}^{n-1} \delta(G_1^i(x)) 
\in H_1(S, R)
\]
if this limit exists.  
\end{defn}

If we are discussing the case when $S$ is an annulus $A$ we will often
identify $G_t$ with a lift $\ti G$ of $G_1$ to the universal cover of
$A$ and identify $H_1(A,\R)$ with $\R.$ We will then refer to the
resulting value of $\rho(x, \ti G)$ as the {\em rotation number} of
$x.$

Even for a general $S$, it follows from the ergodic theorem that the
limit $\rho(x, G_t)$ exists for almost all $x$ with respect to any
invariant Borel measure and defines an integrable function.  It is also clear
that $\rho(x, G_t)$ is independent of the choice of the displacement
function $\delta$ used to define it. This is because the two homology
classes formed by closing the path from $x$ to $G_1^n(x)$ in two
different ways would differ by an amount bounded independently of $n$
and hence this difference disappears in the limit defining $\rho(x,
G_t)$.

We also observe that  the ergodic theorem implies that for any
invariant Borel measure $\nu$ 
\[
\int_S \rho(x, G_t)\ d\nu = \int_S \delta(x)\ d\nu.
\]
This quantity is a key element of our discussion, in case the measure
$\nu$ is the probability measure $\mu$ determined by $\omega$.  More
precisely, we let $\mu$ be the measure on $S$ determined by the volume
form
\[
\frac{1}{\area(S)}\ \omega,
\]
i.e., normalized so $\mu(S) = 1.$

\begin{defn}\label{def: mean rotation vec}
We define the {\em mean rotation vector} of the isotopy $G_t$ to be
\[
\rho_\mu(G_t)  = \int_S \rho(x, G_t)\ d\mu = \int_S \delta \ d\mu 
\in H_1(S, R).
\]
\end{defn}

In \cite{Fr5} it is shown that the mean rotation vector $\rho_\mu \co 
\G_\omega(S) \to H_1(S,R)$ is a homomorphism.  If $S$ is an annulus
$A$, as with the pointwise rotation vector, we will often identify
$G_t$ with a lift $\ti G$ of $G_1$ to the universal cover of $A$ and
identify $H_1(A,\R)$ with $\R.$ We will refer to the resulting value
of $\rho_\mu(\ti G)$ as the {\em mean rotation number}.

\begin{defn}\label{def: Hamiltonian}
An area preserving diffeomorphism $G\co  S \to S$ is called {\em
Hamiltonian} provided there is an isotopy $G_t$ of $S$ with $G_1 = G$
and $G_0 = id$ and $\F([G_t]) = 0 \in H^1(M,R).$ Equivalently (see
Proposition~\ref{prop: mean = flux} below), $G$ is Hamiltonian if there
is an isotopy $G_t$ with mean rotation vector $\rho_\mu (G_t) = 0 \in
H_1(S,R).$ We will denote the group of Hamiltonian diffeomorphisms by
$\H(S)$ or simply $\H.$
\end{defn}

Note that because of Remark~\ref{remark: flux area} it makes no
difference in this definition whether the diffeomorphisms $G_t$ are
area preserving for $t \in (0,1)$.  Also note that any orientation
preserving area preserving diffeomorphism $G \co  S^2 \to S^2$ is
Hamiltonian since any diffeomorphism is isotopic to the identity and
$H^1(S^2,R) = 0.$

The mean rotation vector and the flux homomorphism are closely
related.  In fact their values on an isotopy $G_t$ are essentially
Poincar\'e duals. 

\begin{prop}\label{prop: mean = flux}
Suppose $G_t \in \G_\omega(S)$ and
$\rho_\mu(G_t)\in H_1(S,R)$ is its mean rotation vector.  Then if
$u \in H_1(S,R)$,
\[
\rho_\mu(G_t)\wedge u = \frac{\F(G_t)(u)}{\area(S)}
\]
where $\wedge$ denotes the intersection pairing in $H_1(S,R).$
\end{prop}

\proof
It suffices to prove the result for $u = [\alpha]$ where $\alpha$ is
any simple closed curve in $S$ which is not null homologous, since
these classes generate $H_1(S,\R).$ In particular we may assume the
genus of $S$ is at least one.  Hence we fix such a simple closed curve
$\alpha$ and let $u = [\alpha]$.  Since the intersection pairing with
$u$ is a linear function, it is clear from the definition of rotation
vector that if $\delta$ is a homological displacement function
associated with $[G_t]$ then
\begin{equation}\label{eqn: wedge}
\rho_\mu(G_t)\wedge u = \int \rho(x, G_t)\wedge u\ d\mu 
= \int \delta \wedge u\ d\mu.
\end{equation}
The intersection pairing
with $u$ determines a map $H_1(S,\Z) \to \Z$.  This
in turn defines an infinite cyclic cover $p\co  \bar S \to S$ 
with fundamental group the kernel of the homomorphism
$\pi_1(S) \to H_1(S,\Z) \to \Z$.

Geometrically, this cover is obtained by cutting $S$ along $\alpha$ to
form a surface $S'$ with two boundary components (each a copy of
$\alpha$) and letting $\bar S = S' \times \Z/\sim,$ where $\sim$ is
the relation identifying points of the ``positive'' copy of $\alpha$
in $(S', n)$ with the ``negative'' copy of $\alpha$ in $(S', n+1).$
The generator of the covering translations is $T\co  \bar S \to \bar S$
given by $T(x,n) = (x, n+1).$ Let $S_n$ denote the closure in $\bar S$
of the set $(S',n)/\sim.$ Then each $S_n$ is diffeomorphic to $S'$ and
$T^i(\bar S_n) = S_{n+i}.$ Let $\bar G_t\co  \bar S \to \bar S$ denote
the lift of $G_t$ with $\bar G_0 = id.$

We will let $\alpha_0 = S_{-1} \cap S_0$ and $\alpha_i = T^i(\alpha_0)$
so $\partial S_n = \alpha_{n} \cup \alpha_{n+1}$.  We would like to have
the property that $\bar G_1(S_0) \subset \Int(S_{-1} \cup S_0 \cup S_1)$.    This may not be the case, but it will be true if
$S$ is replaced by a finite cover, specifically by $S(k) = \bar S/ T^k$
for some large $k$, since the infinite cyclic cover $\bar S(k)$ is just
$\bar S$ with covering translations generated by $T^k.$

But both sides of the equality in the statement of the proposition change by a factor of
$1/k$ if we change from $S$ to $S(k)$.  This is because the flux
homomorphism applied to a lift of $[\alpha]$ will not change (see
Remark~\ref{remark: flux cover}) while $\area(S(k)) = k \area(S).$ For
the other side of the equality, the intersection
number $\rho(x, G_t)\wedge [\alpha]$ will decrease by a factor of
$1/k$ for a lift of $\alpha$ in the cover $S(k)$ so
$\rho_\mu(G_t)\wedge [\alpha]$ will decrease by a factor of $1/k$.
Hence the result we want for $S$ is true if and only if it is true for
$S(k)$.  This means we may assume without loss of generality that
$\bar G_1(S_0) \subset \Int(S_{-1} \cup S_0 \cup S_1)$ and hence that  $\bar G_1(S_n) \subset \Int(S_{n-1} \cup S_n \cup S_{n+1})$ for all $n$.   Note in particular that $\bar G_1(\alpha_0) \cap \alpha_i = \emptyset$ for all $i \ne 0$.

To define the homological displacement function for $G_t$ and hence
the rotation vector we considered $S$ as a polygon with edge identifications.
We wish now to choose this data in such a way that the closed curve
$\alpha$ is the image of two identified edges of this polygon.  Then 
the corresponding homological displacement function $\delta$ 
has a lift $\bar \delta\co  \bar S \to H_1(S,R)$ which satisfies
\[
\bar \delta(x) \wedge u = \bar \delta(x) \wedge [\alpha] =
\begin{cases}
	-1, &\text{if $x \in S_n$ and $G_1(x) \in S_{n-1}$;} \\
	0, &\text{if $x \in S_n$ and $G_1(x) \in S_n$;} \\
	+1, &\text{if $x \in S_n$ and $G_1(x) \in S_{n+1}$.}
\end{cases}
\]
It follows that if $D_n^+ = \bar G_1(S_n) \cap S_{n+1}$ 
and $D_n^- = \bar G_1(S_n) \cap S_{n-1}$, then
\begin{align*}
\rho_\mu(G_t)\wedge u &= \int_S \delta \wedge u\ d\mu 
= \int_{S_0} \bar \delta \wedge u\ d\mu \\
&= \mu(D_0^+) -\mu(D_0^-)
= \mu(D_0^+) -\mu(D_1^-)\\
&= \frac{\area(D_0^+) -\area(D_1^-)}{\area(S)}.
\end{align*}
Let $V$ denote the compact submanifold of 
$\bar S$ bounded by $\alpha_0$ and $\bar G_1(\alpha_1)$.  Equivalently
$V = (\bar G_1(S_0) \cap S_1) \cup (S_0 \setminus \Int(G_1(S_1))).$
Then 
\begin{equation}\label{eqn: area diff}
\int_V \bar \omega - \int_{S_0} \bar \omega = \area(D_0^+) -\area(D_1^-)
= \area(S) \big (\rho_\mu(G_t)\wedge u \big ).
\end{equation}
Since $\bar S$ has no de Rham cohomology in dimension $2$, we know that
there is a one-form $\theta$ on $\bar S$ such that $d\theta = \bar \omega$
where $\bar \omega$ is the lift of $\omega.$  We can apply Stokes' theorem
to 
\begin{equation}\label{eqn: theta diff}
\int_V \bar \omega - \int_{S_0} \bar \omega
= \int_{\partial V} \theta\  - \int_{\partial S_0} \theta 
= \int_{\bar G_1(\alpha_1)}\theta \ - \int_{\alpha_1} \theta.
\end{equation}
Now let $G\co  \alpha \times [0,1] \to S$ be defined by $G(x,t) = G_t(x)$
so that by definition 
\[
\F(G_t)(u) = \int_{\alpha \times [0,1]} G^*(\omega).
\]
If $\bar G\co  \alpha \times [0,1] \to \bar S$ is the lift of $G$
with $\bar G(\alpha \times \{0\}) = \alpha_1$ then again applying Stokes'
theorem we have
\[
\int_{\alpha \times [0,1]}\! G^*(\omega) 
= \int_{\alpha \times [0,1]}\! \bar G^*(\bar \omega)
= \int_{\alpha \times \{1\}}\! \bar G^*(\theta)
- \int_{\alpha \times \{0\}}\! \bar G^*(\theta)
= \int_{\bar G_1(\alpha_1)}\!\theta \ - \int_{\alpha_1} \theta.
\]
Combining this with equations~(\ref{eqn: theta diff}) 
and (\ref{eqn: area diff}) we conclude that
$$
\frac{\F(G_t)(u)}{\area(S)} = \rho_\mu(G_t) \wedge u.\eqno{\qed}$$

As an immediate corollary of Proposition~\ref{prop: mean = flux}
we have the following.

\begin{cor} \label{two definitions are equal} If $G_t \in \G_\omega(S)$ and
$\rho_\mu(G_t)\in H_1(S,R)$ is its mean rotation vector, then $\F([G_t]) = 0$ if and only if 
$\rho_\mu (G_t) = 0$.
\end{cor}

The following result is due to Conley and Zehnder \cite{CZ}.

\begin{thm}\label{thm: three pts torus}
If $F\co  T^2 \to T^2$ is a Hamiltonian
diffeomorphism then $F$ has
at least three fixed points.  Moreover three points $\{x_1, x_2, x_3\}
\subset \Fix(F)$ may be chosen so that $F$ is isotopic to the
identity $rel \{x_1, x_2, x_3\}$.
\end{thm}

\section{Chain recurrence}

We will make use of several facts about area preserving surface
homeomorphisms. 

\begin{thm}\label{thm: Brouwer}
An orientation preserving homeomorphism $h$ of the disk $D^2$ which
leaves invariant a measure whose support intersects the interior of the
disk must have an interior fixed point.  In particular if $h$ is 
area preserving then it has an interior fixed point.
\end{thm}
\begin{proof}
The Poincar\'e recurrence theorem says that there is a recurrent point in
the interior of $D^2$.  The Brouwer plane translation theorem then asserts
that the restriction of $h$ to the interior of the disk must have a 
fixed point.
\end{proof}

 Suppose that $h \co  A \to A$ is a homeomorphism of the closed annulus
that is isotopic to the identity.  Recall that a sequence of points
$x_1,\dots, x_n$ in $A$ is an {\em $\varepsilon$-chain from $x_1$ to
$x_n$} if the distance from $h(x_i)$ to $x_{i+1}$ is less than
$\varepsilon$ for $1 \le i \le n-1$.  A point $x$ is {\em chain
recurrent} if there is an $\varepsilon$-chain from $x$ to itself for
all $\varepsilon > 0$.  The set of chain recurrent points is denoted
$\cR(h)$.  We say that $x,y \in \cR(h)$ are in the same {\em chain
transitive component} if for every $\varepsilon >0$ there is an
$\varepsilon$-chain from $x$ to $y$ and another from $y$ to $x$.  If
$\cR(h) = A$ and if all points in $A$ belong to the same chain
transitive component then {\em $h$ is chain recurrent}.

Let $\ti h \co  \ti A \to \ti A$ be a lift to the universal covering
space.  Identify $\ti A$ with $ \mathbb R \times [0,1]$ and denote
projection onto the first coordinate by $p_1$.
The {\em rotation number} of $x$ under $\ti h$ was defined
above to be the number obtained from the rotation vector by
identifying $H_1(A,\R)$ with $\R$.  Equivalently it is
\[
\rho(x,\ti h) = \lim_{n\to\infty} \frac{p_1(\ti F^n(x))-p_1(x)}{n},
\]
if this limit exists.

The following result is an immediate consequence 
of Theorem 2.2 of \cite{Fr1}.

\begin{thm}\label{thm: Birkhoff}
Suppose $h\co  A \to A$ is a chain recurrent homeomorphism of the closed
annulus that preserves orientation and boundary components.  Let $\ti
h \co  \ti A \to \ti A$ be a lift to the universal covering space.  If
$w_0, w_1 \in \ti A$ and $p/q \in \Q$ satisfy
\[
\rho( w_0, \ti h) \le p/q \le \rho( w_1, \ti h)
\]
then $h$ has a periodic point with rotation number $p/q$.  If
$p$ and $q$ are relatively prime this point has period $q$.
\end{thm}

We will use the following method for proving that homeomorphisms are
chain recurrent.

\begin{defn}  A homeomorphism $h \co  A \to A$ of a closed  annulus satisfies  the
{\it essential circle intersection property} provided that every essential
simple closed curve $\beta \subset S$ satisfies $h(\beta) \cap \beta \ne
\emptyset$.
\end{defn}

\begin{lemma} \label{lem: ecip}  Suppose that $A$ is a closed annulus and
$h \co  A \to A$ is a homeomorphism.  Suppose there does not exist a
compact disk $D \subset \Int(A)$ and an $n \ne 0$ such that $h^n(D)
\subset \Int(D)$.  If $h$ satisfies the essential circle intersection
property and if $\partial A$ is in the chain recurrent set of $h$ then
the chain recurrent set of $h$ is all of $A.$
\end{lemma}

\begin{proof}  The homeomorphism $h$ possesses a {\em complete
Lyapunov function} $\phi$: $A \to \R$ (see \cite{ Fr0}).  Recall this is a
smooth function with $\phi(\cR(h))$  compact and nowhere dense 
and satisfying $\phi(h(x)) < \phi(x)$ for $x \notin \cR(h)$.  
In particular $\phi$ is constant if and only if all of $A$ is chain
recurrent.  If $\phi$ is not constant then there is a regular value
$c \in \phi(A)$ with $c \notin \phi(\cR(h))$.  Hence $\phi^{-1}(c)$ is a finite
set of circles in $\Int(D)$, disjoint from their $h$-image.  By hypothesis, these
circles must be inessential and so must bound disks $\{D_j\}$ in
$\Int(A)$. There is no loss in assuming that these disks are components of
$\phi^{-1} (-\infty, c]$.  But then $h(\cup D_j)$ is a proper subset of
$\cup D_j$ and hence for some $n > 0$ one of the $D$'s is mapped into
its interior by $h^n$, a contradiction.
\end{proof}

\section{Hyperbolic structures}

Some of our proofs rely on mapping class group techniques that use
hyperbolic geometry.  In this section we establish notation and recall
standard results about hyperbolic structures on surfaces.  Let $S$ be
a closed orientable surface.  We will say that a connected open subset
$M$ of $S$ has {\em negative Euler characteristic} if $H_1(M, \R)$ is
infinite dimensional or if $M$ is of finite type and the usual
definition of Euler characteristic has a negative value.  Suppose that
$F\co  S \to S$ is an orientation preserving homeomorphism, that $K
\subset S$ is a closed $F$-invariant zero dimensional set and that $S
\setminus K$ has negative Euler characteristic.  Let $M = S \setminus
K$ and $f =F|_M \co M\to M$.

If $K$ is infinite, the surface $M$ can be written as an increasing union of finitely
punctured compact connected subsurfaces $M_i$ whose boundary components
determine essential non-peripheral homotopy classes in $M$.  We may
assume that boundary curves in $M_{i+1}$ are not parallel to boundary
curves in $M_i$.  It is straightforward (see \cite{casble:nielsen}) to put compatible
hyperbolic structures on the $M_i$'s whose union defines a complete
hyperbolic structure on $M$.  Of course, when $K$ is finite, $M$ also has a complete hyperbolic structure.  All hyperbolic structures in this paper are assumed to be complete.

We use the Poincar\'e disk model for the hyperbolic plane $H$.  In
this model, $H$ is identified with the interior of the unit disk and
geodesics are segments of Euclidean circles and straight lines that
meet the boundary in right angles. A choice of hyperbolic structure on
$M$ provides an identification of the universal cover $\ti M$ of $M$
with $H$.  Under this identification covering translations become
isometries of $H$ and geodesics in $M$ lift to geodesics in $H$.  The
compactification of the interior of the unit disk by the unit circle
induces a compactification of $H$ by the \lq circle at infinity\rq\
$\sinfty$.  Geodesics in $H$ have unique endpoints on $\sinfty$.
Conversely, any pair of distinct points on $\sinfty$ are the endpoints
of a unique geodesic. 

 Each covering translation $T$ extends to a homeomorphism (also called) $T \co  H \cup \sinfty \to H \cup \sinfty$. The fixed point set of a non-trivial $T$ is either one or two points in $\sinfty$. We denote these point(s) by $T^+$ and $T^-$, allowing the possibility that $T^+ = T^-$.  If $T^+ = T^-$, then $T$ is said to be {\it parabolic}.  If $T^+$ and $T^-$ are distinct, then    $T$ is said to be {\it hyperbolic} and we may assume that $T^+$ is a sink and $T^-$ is a source.   

We use the identification of $H$ with $\ti M$ and write $\ti f \co  H \to
H$ for lifts of $f\co  M \to M$ to the universal cover.  A fundamental
result of Nielsen theory is that every lift $\ti f \co  H \to H$ extends
uniquely to a homeomorphism (also called) $\ti f \co  H \cup \sinfty \to
H \cup \sinfty$.  (A proof of this fact appears in Proposition 3.1 of
\cite{han:fpt}).  If $f \co  M \to M$ is isotopic to the identity then
there is a unique lift $\ti f$, called the {\it identity lift,} that
commutes with all covering translations and whose extension over
$\sinfty$ is the identity.

 For any extended lift $\ti f \co  H \cup \sinfty \to
H \cup \sinfty$ there is an {\it associated action $\ti f_\#$ on
geodesics in $H$ }defined by sending the geodesic with endpoints $P$
and $Q$ to the geodesic with endpoints $\ti f(P)$ and $\ti f(Q)$.
The action $\ti f_\#$ projects to an {\it action $f_\#$ on geodesics
in $M$}.

 We occasionally allow $\partial S$ to be non-empty.
There is still a hyperbolic structure on $M$ but now the universal
cover $\ti M$ is naturally identified with the intersection of $H$ with the interior of the convex
hull of a Cantor set $C \subset \sinfty$. The frontier of $\ti M$ in $H \cup \sinfty$ is
the union of $C$ with the full pre-image of $\partial S$.  If $f$ is
isotopic to the identity, then the identity lift commutes with all
covering translations and extends to a homeomorphism of the frontier
of $\ti M$ that fixes $C$.

\section{Lifting to annuli}

Throughout this section $F\co  S \to S$ is a diffeomorphism of a closed
surface $S$.  An isolated end of an open set $U \subset S$ has
neighborhoods of the form $N(E) = S^1 \times [0,1)$.  The set $fr(E) =
cl_S(N(E)) \setminus N(E)$, called the {\it frontier of $E$}, is 
independent of the choice of $N(E)$.

\begin{lemma}\label{lem: extension}   
Suppose that $M \subset S$ is an open $F$-invariant set and that $E$ is an
isolated end of $M$ whose frontier is contained in Fix($F$).  Then $E$
can be compactified by a circle $C$ and $f = F|_M$ can be extended
continuously over $C.$ If the frontier of $E$ is not a single point
then $f|_C =$ identity.
\end{lemma}

\begin{proof} 
If $fr(E)$ is a single point of Fix($F$) then we compactify by blowing
up that point using the
differentiability of $f$.  We may therefore assume that $fr(E)$ is a non-trivial
compact connected subset of $\Fix(F)$.  The existence of
$C$ and a continuous extension of $f$ over $C$ is a consequence of
the theory of prime ends (see \cite{M2} for a good modern exposition).
Moreover, a dense set  $P \subset C$ corresponds to accessible points
in $fr(E)$, i.e.\ for any $p \in P$ there is an embedded
 arc $\alpha\co  [0,1] \to M \cup C$ with $\alpha(0) = p$ and
$\alpha((0,1]) \subset M$ and such that the restriction of $\alpha$ to
$(0,1]$ extends to a continuous map $\hat \alpha\co  [0,1] \to S$ with
$\hat \alpha(0) \in fr(E).$ A result of Mather (Theorem 18 of
\cite{M2}) asserts that two such extensions, $\hat \alpha_1$ and $\hat
\alpha_2$ correspond to the same point of $C$ if and only if $\hat
\alpha_1(0) = \hat \alpha_2(0)$ and $\hat \alpha_1$ and $\hat
\alpha_2$ are isotopic by an isotopy with support in $M.$ A
consequence of this is that the homeomorphism $f|_C$ depends only on
the isotopy class of $F$ relative to $fr(E)$.  By Lemma 4.1 of
\cite{han:commuting} $F$ is isotopic rel $fr(E)$ to a homeomorphism
that is the identity on a neighborhood of $fr(E)$.  Thus $f|_C$ is the
identity as desired.
\end{proof}

\begin{lemma} \label{canonical cyclic lift}
Suppose that $M$ is a connected $F$-invariant open subset of $S$ whose
frontier, $fr(M) = cl(M) \setminus M$, is contained in Fix$(F)$.  If
$S=T^2$, assume that $fr(M)\ne \emptyset$.  Denote $F|_M$ by $f$.  If
$\beta \subset M$ is a simple closed curve that is essential in $M$
and if $f(\beta)$ is homotopic to $\beta$ in $M$ then there is a
closed annulus $\A$, a homeomorphism $\hat f \co \A \to \A$ that is
isotopic to the identity and a covering map $\pi\co  \Int(\A) \to M$ such
that
\begin{enumerate}
\item $\hat f|_{\Int \A}$ is a lift of $f$.
\item $\pi^{-1}(\beta)$ contains an essential simple 
closed curve $\hat \beta$.  
\item if $\beta$ is non-peripheral    then  $\hat f |_{  \partial \A}$ depends only on the isotopy class of $f$. In particular, if $\beta$ is non-peripheral and $f\co  M \to M$ is isotopic to the identity  then $\hat f |_{  \partial \A}$ is  the identity. 
\item if $f\co  M \to M$ is isotopic to the identity,  and   either $S \ne S^2$ or $S \setminus M$ contains at least three points  then  $\hat f |_{  \partial \A}$  has fixed points 
\end{enumerate}
We will call $\pi$ the {\em canonical cyclic cover} associated to
$\beta$ and $\hat f$ the {\em canonical cyclic lift} of $f$.
\end{lemma}

\begin{proof} 
If $M$ is an open annulus then $\Int(\A) = M$, $\hat f|_{\Int(\A)} = f$
and it suffices to compactify the ends of $M$ by circles and extend
$f$ over these circles.  Lemma~\ref{lem: extension} therefore
completes the proof in this case.  We may now assume that $M$ has 
negative Euler characteristic.

Suppose at first that $\beta$ is not peripheral.  Let $H$ be the universal cover of $M$,  let $\ti \beta$ be a lift of $\beta$ and let $T\co  H \cup
\sinfty \to H \cup \sinfty$ be the indivisible covering translation that fixes the endpoints of $\ti \beta$.
Since $f(\beta)$ is homotopic to $\beta$, there is an extended lift $\ti f \co  H \cup
\sinfty \to H \cup \sinfty$  that  commutes with $T$.  The choice of $\ti f$ is well defined up to composition with iterates of $T$ because these are the only covering translations that commute with $T$.  Let $\A$ denote the quotient space obtained by dividing $(H \cup S_{\infty}) \setminus T^{\pm}$ by
the action of $T$.  Denote by $\hat f\co  \A \to \A$ the homeomorphism
induced by $\ti f$.  It is independent of the choice of lift $\ti f$,
so $\hat f$ is unique. 
The extension $\ti f|_{S_{\infty}}$ depends only on the isotopy class of
$f$ and the choice of lift $\ti f$.  Since $\hat f$ is independent of
the choice of lift, $\hat f |_{\partial \A}$ depends only on the isotopy
class of $f$.  If $f$ is isotopic to the identity then we may take
$\ti f$ to be the identity lift.  Thus $\ti f |_{\sinfty}$ and $\hat
f|_{\partial \A}$ are the identity.

If $\beta$ is peripheral, compactify the isolated end corresponding to
$\beta$ using Lemma \ref{lem: extension} and repeat the previous
argument.  If $f$ is isotopic to the identity then the identity lift
fixes a Cantor set in $\sinfty$ so $\hat f$ fixes a Cantor set in
$\partial \A$.
\end{proof}

The following lemma relates periodic points of $\hat f$ to periodic points of $f$.

\begin{lemma} \label{not on the boundary}
Suppose that $\hat f \co  \A \to \A$ is the canonical cyclic lift of $f \co 
M \to M$ associated to a simple closed curve $\beta$ and that $\hat x$
is a periodic point of $\hat f$ with period $q$.  If $\rho(\hat x,
\hat f)$ is non-zero and does not occur as $\rho(\hat z, \hat f)$ for
some $\hat z \in \partial \A$, then the image $x \in M$ of $\hat x$ is
periodic for $f$ with period $q$.
\end{lemma}

\begin{proof}  We use the notation of the proof of Lemma~\ref{canonical cyclic lift}.  If $M$ is an open annulus then the lemma is obvious so   assume that $M$ has negative Euler characteristic.  There is a lift $\ti x \in H$ of $ \hat x$  such that $\ti f^q(\ti x) = T^p(\ti x)$ where $\rho(\ti x, \ti f) =p/q$.  Since $f^q(x) = x$, it suffices to show that $  f^s(x) \ne x$ for  $1 \le s \le q-1$.  Suppose to the contrary that   there exists $1 \le s \le q-1$ and a covering translation $S \co  H \to H$ such that $\ti f ^s(\ti x) = S(\ti x)$.  Recall that the covering translation $\ti f^q S \ti f^{-q}$  is denoted  $\ti f^q_\#(S)$.  We have  
\begin{align*}
\ti f^q_\#(S) (T^{p}(\ti x)) &= \ti f^q_\#(S) (\ti f^{q}(\ti x)) = \ti
f^{q}(S(\ti x)) = \ti f^{q + s}(\ti x) \\ &= \ti f^s \ti f^{q}(\ti x)
= \ti f^s T^{p}(\ti x) =T^{p}(\ti f^s (\ti x)) = T^{p}(S(\ti x)).
\end{align*}
Covering translation that agree on a point are equal so $\ti f^q_\#(S)
= T^p S T^{-p}$.  Let $z \in S_{\infty}$ be the unique attracting
fixed point for the action of $S$.  Then $T^p(z)$ and $\ti f^q(z)$ are
both the unique attracting fixed point for the action of $\ti
f^q_\#(S) = T^p S T^{-p}$.  Thus $\ti f^q(z) = T^p(z)$.  We are
assuming that $\hat f^s(\hat x) \ne \hat x$ and hence that $S \ne T^j$
for any $j$.  In particular, $z \ne T^{\pm}$ and $z$ projects to a
point $\hat z \in \partial \A$ such that $\rho(\hat z, \hat f) = p/q$.
This contradiction completes the proof.
\end{proof}

We use the following lemma in several places to prove the existence of periodic points of arbitrarily high period.

\begin{lemma}  \label{lem: per points f hat} Suppose that
$F \in \H(S)$ and that $M$ is a connected $F$-invariant open subset of $S$ whose
frontier, $fr(M) = cl(M) \setminus M$, is contained in Fix$(F)$.  If
$S=T^2$, assume that $fr(M)\ne \emptyset$.  Denote $F|_M$ by $f$.  Suppose also that
$\beta \subset M$ is a simple closed curve whose homotopy class in $M$ is
non-trivial and fixed by $f$ and that the lift 
$\hat f \co  \A \to \A$ to the canonical
cyclic cover associated to $\beta$  has points with 
two distinct rotation numbers.  Then there are periodic points of $ f$
of every sufficiently high period.  
\end{lemma}

\begin{proof}   Let $\pi \co  {\rm int}(\A) \to M$ be the covering map
and let $\hat \omega = \pi^*(\omega)$.  
The smooth measure on the interior of $\A$ determined  by $\hat \omega$ is $\hat f$-invariant. 

We want first to show that $\hat f \co  \A \to \A$    has the essential circle
intersection property.  Suppose that $\hat \alpha$ is an essential simple
closed curve in $\A$ and that $\hat f(\hat \alpha) \cap \hat  \alpha =
\emptyset$.  There is no loss in assuming  that $\hat \alpha$ is smooth
and contained in the interior of $\A$.   

The region bounded by $\hat
\alpha$ and $\hat f(\hat \alpha)$ is an annulus that we denote $ \hat
A_0$. Choose a smooth parametrization  $\hat h\co  S^1 \times [0,1] \to \hat A_0$
such that $\hat h(s,1) = \hat f(\hat h(s,0))$ for each $s \in S^1$ and
let $h\co  S^1 \times [0,1] \to M$ be $h = \pi \circ \hat h.$
Then
\begin{equation}\label{eqn: h omega}
\int_{S^1 \times [0,1]} h^*(\omega)
= \int_{S^1 \times [0,1]} \hat h^*(\hat \omega)
= \int_{\hat A_0} \hat \omega \ne 0.
\end{equation}
On the other hand since $F \in \H(S)$ there is an isotopy $F_t$ with
$F_0 = id$ and $F_1 = F.$ In the case $S = T^2$ or $S^2$ the frontier
of $M$ is non-empty and consists of fixed points.  Hence we may assume
there is a point $p \in fr(M)$ such that $F_t(p) = p$ for all $t,$ since any
diffeomorphism of $S^2$ or $T^2$, which is isotopic to the identity and has a
fixed point $p,$ is isotopic to the identity relative to $p$.

Let $\alpha = \pi \circ \hat \alpha \co  S^1 \to M$. 
If $F_\alpha \co  S^1 \times [0,1] \to S$ is given by 
$F_\alpha (s, t) = F_t(\alpha(s))$ then
\begin{equation}\label{eqn: FF omega}
\int_{S^1 \times [0,1]} F_\alpha^*(\omega)
= \F(F_t)([\alpha]) = 0
\end{equation}
by hypothesis.

The two maps $h, F_\alpha \co  S^1 \times [0,1] \to S$ 
agree on their ends, so we may glue them together to
form a map on the torus, $G\co  T^2 \to S.$  The map $G$
must have degree zero because any map of a torus into a surface
of negative Euler characteristic has degree zero and
in the cases that $S = S^2$ or $S = T^2$ we have the fact that
the image of $G$ lies in $S \setminus \{p\}.$
Consequently
\[
0 = \int_{T^2} G^*(\omega)
= \int_{S^1 \times [0,1]} F_\alpha^*(\omega)
- \int_{S^1 \times [0,1]} h^*(\omega)
\]
which contradicts equations~(\ref{eqn: h omega}) and (\ref{eqn: FF omega}).
Thus we have shown that $\hat f$ has the essential circle intersection
property.

We would like to show that $\hat f$ is chain recurrent.  If we knew
that $\partial \A \subset \cR(\hat f)$ then this would follow from
Lemma~\ref{lem: ecip}.  It turns out we can reduce to the case
$\partial \A \subset \cR(\hat f)$ by extending $\hat f$ to a slightly
larger annulus without affecting the periods of periodic points.  If
the rotation number of the restriction of $\hat f$ to a boundary
component is irrational then that component is in $\cR(\hat f)$ since
this circle homeomorphism must be conjugate to either an irrational
rotation or a Denjoy type example, both of which have all points chain
recurrent.  If the rotation number on a boundary circle is rational
there is a periodic point $p_0$ on that circle.  In this case there is
an isotopy $h_t$ on the circle, with $h_0$ the given homeomorphism
and $h_1$ a finite order homeomorphism, and which has the property
that $\Per(h_t) = \Per(h_0)$ for all $t \in [0,1).$
We now attach a collar neighborhood $S^1 \times [0,1]$ to
the boundary of $A$ and extend $\hat f$ to the union by having it
preserve concentric circles of the collar and act on them in the way
prescribed by the isotopy $h_t$.  We can do this to both boundary components
if necessary.  Call the enlarged annulus $\B$ and the extended map 
$\hat g\co \B \to \B$.  Then $\partial \B \subset \cR(\hat g)$ and $\hat g$ has
a point of period $q$ with rotation number $\frac{p}{q}$ if and only if there is a point of period $q$ with rotation number $\frac{p}{q}$ 
for $\hat f.$

  It is also the case that $\hat g$ satisfies the essential circle
intersection property. This is because if $\gamma$
is an essential circle with $\hat g(\gamma) \cap \gamma = \emptyset$
and $x \in \gamma$ then points $z \in \alpha(x)$ and $w \in \omega(x)$
lie in $\cR(\hat g)$ but in different components of $\cR(\hat g)$
since $\gamma$ separates them.  On the other hand $\cR(\hat g) \cap
(\B \setminus \A)$ consists of all periodic points in $\B \setminus
\A$ and has two components, one containing each of the components of
$\partial \B.$ It follows that $\gamma$ cannot intersect $\B \setminus
\A$ since any point of $\B \setminus \A$ has both its alpha and omega
limit sets in the same component of $\cR(\hat g) \cap (\B \setminus
\A).$  Thus $\gamma \subset A$ and since $\hat f$ has the essential circle
intersection property, so does $\hat g.$ Also
if $D$ is a closed disk and $\hat g^n(D) \subset \Int(D)$ then $D \cap
(\B \setminus \A) = \emptyset$ because the omega limit set of any
point in $\B \setminus \A$ lies in the component of $\Per(h_0)$
containing $\partial \B.$

By Lemma~\ref{not on the boundary} it suffices to show that every sufficiently large $q$ occurs as the period of a periodic point for $\hat f$ with non-zero rotation number not realized on $\partial \A$.  This property holds for $\hat f$ if and only if holds for $\hat g$ so  Lemma~\ref{lem: ecip} and Theorem~\ref{thm: Birkhoff} complete the proof.  
\end{proof}

\section{Normal form}\label{normal_form}

     In this section we introduce a normal form for certain mapping classes of infinitely punctured surfaces.

\begin{defn} \label{can form} 
An orientation preserving homeomorphism $F\co  S \to S$ of an orientable
closed surface {\it has a normal form relative to its fixed point set}
if there is a finite set $R$ of simple closed curves in $M = S
\setminus \Fix(F)$ and a homeomorphism $\phi$ isotopic to $F$ rel $\Fix(F)$
such that:
\begin{description}
\item [(1)] $\phi$ permutes disjoint open annulus neighborhoods $A_j
\subset M$ of the elements $\gamma_j \in R $.
\end{description}

Let $\{S_i\}$ be the components of $S \setminus \cup A_j$, let $X_i =
\Fix(F) \cap S_i$, let $M_i = S_i \setminus X_i$ and let $r_i$ be the
smallest positive integer such that $\phi^{r_i}(M_i) = M_i$ .  Note that $r_i=1$ if $X_i \ne \emptyset$.

\begin{description}
\item [(2)] If $X_i$ is infinite then $\phi|_{S_i}$= identity. 
\item [(3)] If $X_i$ is finite then $M_i$ has negative Euler characteristic and $\phi^{r_i}|_{M_i}$  is either pseudo-Anosov  or periodic.  In the periodic case, $\phi^{r_i}|_{M_i}$ is an isometry of a hyperbolic structure on $M_i$.

\end{description}
We say that $\phi$ is a {\it normal form for $F$} and that $R$ is {\it
the set of reducing curves for $\phi$}.
\end{defn}

If  $R$ has minimal cardinality among all sets of reducing curves for all normal forms for $F$, then we say that {\it $R$ is a minimal set of reducing curves}.  The following lemma asserts that, up to isotopy rel Fix($F$), there is a unique minimal set $R(F)$ of reducing curves.

\begin{lemma}  \label{unique} 
If $F$ has a canonical form relative to its fixed point set then the
minimal reducing set $R(F)$ is well defined up to isotopy rel $\Fix(F)$.
\end{lemma}

\begin{proof} Suppose that $R$ is a
minimal reducing set for the normal form $\phi$. Let $Y(F)$ be a
finite subset of Fix($F$) that contains each finite $X_i$ and at least
two points from each infinite $X_i$.  Then $\phi$ is a Thurston normal
form \cite{Th} for the isotopy class of $F$ relative to
$Y(F)$ and $R$ is a minimal set of reducing curves for this relative
isotopy class.  It is shown in section 2 of \cite{han:nielsen} that $R$
is well defined up to isotopy rel $Y(F)$. Since this holds for all
$Y(F)$, $R$ is well defined up to isotopy rel Fix($F$).
\end{proof}

\noindent{\bf Theorem~\ref{canonical form}}\qua {\sl 
Every orientation preserving diffeomorphism $F \co  S \to S$ of a closed
orientable surface has a normal form up to isotopy relative to its
fixed point set.}

\begin{proof}
By Lemma 4.1 of \cite{han:commuting} there is a neighborhood $W$ of the accumulation set $A(F)$ of Fix($F$) such that:
\begin{itemize}
\item $W$ and $S \setminus W$ have  finitely many components.
\item $F|_W \co  W \to S$  is isotopic to the inclusion relative to Fix($F)\cap W$.
\item if $x \in W\setminus \Fix(F)$ then the path from $F(x)$ to $x$
determined by this isotopy is contained in $S\setminus \Fix(F)$. 
\end{itemize}

We may assume without loss that each component of $W$ intersects
Fix($F$) in an infinite set.  The isotopy extension theorem (Theorem 5.8 of \cite{mil}) and the second and third item imply that $F$
is isotopic rel Fix($F$) to $\phi_1 \co S \to S$ satisfying $\phi_1|_W
= {\rm identity}$.  By \cite{brn-kis} each component of $S \setminus W$ is  invariant under the action of $\phi_1$.  
 If some component of $S \setminus W$ is a disk that contains at most one element of Fix($F$) then add it to $W$.  After a
further isotopy rel Fix($F$), we may assume that $\phi_1|_W$ is still
the identity.   Suppose that some component of $S \setminus W$ is an open annulus
 $A$ that is disjoint from Fix($F$).  If $\phi_1|_{(W \cup A)}$ is isotopic to
the identity rel Fix($F$), then add $A$ to $W$.  After a further
isotopy rel Fix($F$), we may assume that $\phi_1|_W$ is still the
identity.  If $\phi_1|_{(W \cup A)}$ is not isotopic to the identity rel
Fix($F$) then after a further isotopy rel Fix($F$) we may assume that
$\phi_1|_A$ is a non-trivial Dehn twist.

  Denote $S \setminus \Fix(F)$ by $M$.  At this point $M \setminus (W
\cap M)$ consists of a finite number of unpunctured annuli on which
$\phi_1$  acts by Dehn twists and
finitely many, finitely punctured subsurfaces of negative Euler
characteristic.  Apply the Thurston decomposition theorem to each of
the non-annular subsurfaces to produce $\phi \co  S \to S$.  Let $R_1$ be the union of the core curves of the Dehn twist annuli and   Let $R$ be the union of $R_1$ with the boundary curves of $W$ that are not boundary curves of a Dehn twist annuli.  (Those curves are already accounted for in $R_1$.)   After modifying $\phi$ on annular neighborhoods of the
reducing curves we may assume that (1) is satisfied.  Properties (2)
and (3) are immediate from the construction.
\end{proof}

We now restrict our attention to the case that $F$ is isotopic to the identity. 

\begin{lemma} \label{s2}  Assume that $\phi$, $R$, $M_i$ and $r_i$ are as in Definition~\ref{can form} and that  $F$ is isotopic to the identity.  Then:
\begin{itemize}
\item each  $\gamma \in R$ is $\phi$-invariant.
\item each $r_i=1$.
\item   if $\phi|_{M_i}$ is periodic then it is the identity. 
\end{itemize} \end{lemma}

\begin{proof}
If $\gamma$ bounds a disk $D$ in $S$ then $D \cap \Fix(F)$ is
non-trivial.  The disk $F(D)$, and hence the disk $D'$ bounded by
$\phi(\gamma)$, has the same intersection with $\Fix(F)$ as does $D$.
Since $\phi(\gamma)$ and $\gamma$ are either equal or disjoint and not
parallel in $M$, they must be equal. This completes the proof of the
first item if $\gamma$ is inessential in $S$.  Since $D$ contains
fixed points, $\phi$ does not interchange the components of $S
\setminus \gamma$.  In particular, if $M_i$ is adjacent to the annulus
$A_j$ containing $\gamma$, then $M_i$ is $\phi$-invariant.  This
proves the second item if some  component of $\partial M_i$ is inessential in
$S$.

Suppose that $M_i$ is $\phi$-invariant and that $X$ and $Y$ are distinct
components of $\partial M_i$ that are essential in $S$.  If $Y =
\phi(X)$ as sets then $X$ and $Y$ must be parallel in $S$. Up to
isotopy in $S$, $\phi$ reverses the orientation on the simple closed
curve $X$, in contradiction to the assumption that $F$ is isotopic to
the identity.  We have now shown that if $M_i$ is $\phi$-invariant
then every component of $\partial M_i$ is $\phi$-invariant.

The lemma is obvious if Fix($F) = \emptyset$ so assume that $M$ has at
least one puncture.  Each puncture in $M$ is contained in exactly one
$M_i$ which must therefore be $\phi$-invariant.  If $M_i$ is
$\phi$-invariant then each annulus $A_j$ adjacent to $M_i$ is
$\phi$-invariant and so each $M_k$ adjacent to $A_j$ is
$\phi$-invariant.  The first two items follow by induction.

For the third item we may assume that $M_i$ is finitely punctured. For this argument it is convenient to replace the punctures in $M_i$ with boundary components to form a compact surface $\hat M_i$.  If
$\phi|_{M_i}$ is periodic, then $\phi$ extends to an isometry of a hyperbolic
structure on $\hat M_i$ that setwise fixes each boundary component.  Since
$F$ induces the identity on the first homology of $S$, $\phi|_{\hat M_i}$
induces the identity on first homology of $\hat M_i$.  Any isolated fixed
point of $\phi|_{\hat M_i}$ has positive fixed point index so the Lefschetz
theorem implies that not all fixed points are isolated.  Since
$\phi|_{\hat M_i}$ is an orientation preserving isometry, it must be the
identity.
\end{proof}

\begin{remark} \label{DT} If each $r_i = 1$ and if each $\phi|_{M_i}$ is the identity, then we may assume that $\phi|_ {A_j}$ is an iterate of a Dehn twist about $\gamma_j$.
\end{remark}

\section{A special case}

In this section we prove a special case of Theorem~\ref{thm: periodic point}.

\begin{prop}  \label{prop: per point} Suppose $F\co  S \to S$ is a non-trivial, Hamiltonian
diffeomorphism of a closed oriented surface $S$ and that if $S = S^2$ then 
$\Fix(F)$ contains at least three
points.  If $F$  is
not isotopic to the identity relative to $\Fix(F)$, then there
exist $p>0$ such that $F$ has a periodic point of
 period $k$ for every $k\ge p.$
\end{prop}

\begin{proof}
Let $\phi$ be a canonical form for $F$ with minimal reducing set
$R(F)$.  Let $Y$ be a finite subset of Fix($F$) that contains all
finite $X_i$ and at least two points from each infinite $X_i$. Then
$\phi$ is a Thurston normal form for $F$ relative to $Y$. If some
$\phi|_{M_i}$ is pseudo-Anosov then (\cite{Th}, Expose 10 of \cite{FLP}) there exist $p >0$
such that $F$ has a
periodic point of period $k$ for every $k\ge p$ and we are done. 

By Lemma~\ref{s2} and Remark~\ref{DT}, $\phi$ is a composition of
non-trivial Dehn twists along the elements of $R(F)$.  Thus it
suffices to show that if $R(F)$ is non-empty, then there exist $p>0$ such that $F$ has a periodic point of
 period $k$ for every $k\ge p.$  So
suppose $\beta \in R(F)$.  We work with $M_Y = S \setminus Y$, $f_Y =
F|_{M_Y}$ and $\phi_Y = \phi|_{M_Y}$ .  Let $\hat f_Y,\hat \phi_Y \co  \A
\to \A$ be the canonical cyclic lifts corresponding to $\beta \subset
M_Y$.

We wish to apply Lemma~\ref{lem: per points f hat} to $\hat f_Y$ which
will give points of every sufficiently large period for $f_Y$.
It suffices to show that the rotation numbers of lifts of boundary
points for any lift of $\hat f_Y\co  \A \to \A$ to its universal covering
space are different for the two components of the boundary.  By
Lemma~\ref{canonical cyclic lift}(2), $\hat f_{Y|\partial \A} = \hat
\phi_{Y|\partial \A}$ so we can calculate the rotation number using
$\phi$ in place of $F$.  Let $A_j$ be the annulus containing $\beta$
and let $\beta_0$ and $\beta_1$ be the boundary components of $A_j$.
Let $M_0$ and $M_1$ be the (possibly equal) pieces in the Thurston
decomposition whose boundaries contain $\beta_0$ and $\beta_1$
respectively.

 If $\alpha_0 \ne \beta_0$ is a component of  $\partial M_0$ then
we can join $\alpha_0$ and $\beta_0$ by an arc on which $\phi$ is the
identity.  We can lift this arc and $\alpha_0$ and $\beta_0$ to $H$
in such a way that the lift $\ti \beta_0$ has the same endpoints as
the lift $\ti \beta$ of $\beta$.  There is a lift of $\phi$
which is the identity on $\ti \beta_0,\ \ti \alpha_0$ and the lifted
arc joining them.  Hence it is clear that the endpoints of the lift
$\ti \alpha_0$ which lie in $\sinfty$ will have the same rotation
number as any point of $\ti \beta_0$ for this lift of $\phi$ and thus
for any lift of $\phi.$ If there are no components of the boundary of
$M_0$ other than $\beta_0$ then we may choose an arc in $H$ joining a
point of $\ti \beta_0$ to a point in $\sinfty$ which crosses no lift
of a reducing curve.  Again, one lift of $\phi$ will be the identity
on this arc and on $\ti \beta_0.$ So again points of the boundary of
$\sinfty$ which lie on the same side of $\ti \beta$ as $\ti \beta_0$
all have the same rotation number as points of $\ti \beta_0$ for any
lift of $\phi$.

The same argument shows points in $\sinfty$ on the other side of $\ti \beta$
have the same rotation number as points of $\ti \beta_1$.
It follows that points in the two boundary components of the universal
covering of $\A$ have rotation numbers related like the rotation
numbers of opposite sides of an annulus with a non-trivial Dehn
twist.  i.e.\ they differ by a non-zero integer.  The proof now follows from  Theorem~\ref{thm: Birkhoff} and Lemma~\ref{lem: per points f hat}. 
\end{proof}

\section{Translation and Reeb classes}

We will denote the orbit, backward orbit and forward orbit of $x$ by
$\calo(x)$, $\calo^-(x)$ and $\calo^+(x)$ respectively.  Unless stated
otherwise, all maps are assumed to be orientation preserving
homeomorphisms.

  In this section we recall and enhance some results from \cite{han:fpt}. 

\begin{defn}  Let $Tr \co  \rtwo \to \rtwo$ be translation to the right by one unit.  If $f \co  \rtwo \to \rtwo$ is an orientation preserving homeomorphism and if $X = \cup_{i=1}^r\calo(x_i)$ for $x_i \in \rtwo$ with distinct orbits, then we say  that {\it $f$ is a translation class relative to $X$} if the isotopy class of $f$ relative to $X$ is conjugate to the isotopy class of $Tr$ relative to the union of some (and hence any - see Lemma 2.1 of \cite{han:fpt}) $r$ distinct orbits of $Tr$.  
\end{defn}

\begin{defn}  Let $R \co  \rtwo \to \rtwo$ be a homeomorphism that preserves $\rone \times \{0,1\}$ and that agrees with translation  to the left [respectively right] by one unit on $\rone \times \{1\}$ [respectively $\rone \times \{0\}]$.  If $f \co  \rtwo \to \rtwo$ is an orientation preserving homeomorphism and if $X = \calo(x) \cup \calo(y)$ for some $x,y \in \rtwo$ with distinct $f$ orbits then we say that {\it $f$ is a Reeb class relative to $X$}  if the isotopy class of $f$ relative to $X$ is conjugate to the isotopy class of $R$ relative to the union of the $R$-orbits of $(0,0)$ and $(0,1)$.
\end{defn}

     Parts (1) and (2) of the following theorem are Corollary 2.2 and Theorem 2.6 of  \cite{han:fpt} respectively.

\begin{thm}  \label{fpt} If   $f \co  \rtwo \to \rtwo$ is fixed point free then
\begin{description} 
\item [(1)]  $f$ is a translation class relative to any one of its orbits.
\item [(2)]  $f$ is either a translation class or a Reeb class relative to any two of its orbits.
\end{description}
\end{thm}

      Several of our arguments proceed by contradiction and so require a method for proving that the conclusions of Theorem~\ref{fpt} fail.  Lemma~\ref{intersection} below provides this method.

\begin{defn}  Suppose that $S = \rtwo$ or $S = {\rm int}({\mathbb A})$ and that $f \co  S \to S$ is an orientation preserving homeomorphism.  An arc $\beta' \subset S$ connecting $x $ to $f(x)$ is called a {\it translation arc for $x$} if $f(\beta') \cap \beta' = f(x)$. For our applications there will be a closed infinite discrete set $X \subset S$ that is the union of orbits and that intersects $\beta'$ exactly in $x$ and $f(x)$.  Equip $M = S \setminus X$  with a complete hyperbolic structure.  We say that the unique geodesic $\beta \subset M$ in the homotopy class of $\beta'$ is a {\it \tag}\ for $x$ relative to $X$.
\end{defn}

\begin{remark}  In \cite{han:fpt} a more general object called a homotopy translation arc is considered but for our present purposes this simpler definition will suffice.
\end{remark}

\begin{defn}An open disk $U$ is called a {\it free disk} for the homeomorphism $f $ if $ f(U) \cap U = \emptyset$.
\end{defn}

\begin{lemma}  \label{intersection} Assume that  $f \co  \rtwo \to \rtwo$ is a fixed point free orientation preserving homeomorphism.  Then:
\begin{description}
\item [(1)] If $X$ is a single orbit and $f$ is a translation class relative  to $X$ then \tag s for elements of $X$ are unique.
\item [(2)]  If $X = \calo(x_1) \cup \calo(x_2)$ and $f$ is a Reeb class relative to  $X$ then \tag s for  $x \in X$ are unique. Moreover:
\begin{itemize}
\item no free disk contains two elements of $X$.
\item a \tag\ for $f^i(x_1)$ can not intersect a \tag\ for $f^j(x_2)$.
\end{itemize}
\item [(3)] If $X = \calo(x_1) \cup \calo(x_2)$ and $f$ is a translation class relative to  $X $ then there is a compact set $K$ with the following property:  if $\beta_1$ is a \tag\  for an element  of $\calo^-(x_1)$ and $\beta_2$ is a \tag\  for an element  of $\calo^+(x_2)$ and if both $\beta_1$ and $\beta_2$ are disjoint from $K$ then $\beta_1 \cap \beta_2 = \emptyset$.
\end{description}
\end{lemma}

\noindent{\bf Proof of Lemma~\ref{intersection}}\qua   (1) is Theorem 2.2   of \cite{han:fpt}. 

To prove (2) we will prove the a priori stronger statement that there is a unique geodesic arc $\alpha \subset M$ connecting $x$ to $f(x)$ with
the property that $f_\#(\alpha) \cap \alpha = f(x)$.  (It is immediate from the definitions that this  property depends only on the conjugacy class of the relative isotopy class of $f$ which is not the case if we insist that $\alpha$ is a \tag.)    It suffices to
assume that $f = R$, $x_1= (0,0)$ and $x_2 = (0,1)$.  We will give the
argument for $x_2$, all other cases being completely analogous.  We
may assume that $\rone \times \{0,1\}$ is a union of geodesics.
Suppose that $\alpha \ne [0,1]\times \{1\}$.  After an isotopy that
preserves $\rone \times \{0,1\}$, we may assume (c.f. Lemma 3.5-(3) in
\cite{han:fpt}) that $R(\alpha)$ is a geodesic.  In particular
$R(\alpha)$ and $\alpha$ intersect only in their common endpoint.  Let
$\alpha_0$ be the maximal initial segment of $\alpha$ whose interior
is disjoint from $\rone \times \{0,1\}$ and let $w$ be the terminal
endpoint of $\alpha_0$.  Since regions bounded by a pair of geodesic
arcs must contain punctures, $\alpha_0 \ne \alpha$.  In particular
$R(\alpha_0)$ is disjoint from $\alpha_0$.  There are three
possibilities: $w \in (-1,1) \times \{1\}$; $w \in (k,k+1) \times
\{1\}$ for some $k \ne -1,0$; and $w \in \rone \times\{0\}$.  The
first option is impossible since the region bounded by a pair of
geodesic arcs must contain punctures.  The latter two are impossible
because in those cases $R(\alpha_0)$ could not be disjoint from
$\alpha_0$. This proves the main statement in (2).  The first bulleted
item follows from the standard construction of translation arcs (see,
for example, Lemma 4.1 in \cite{han:fpt}) since the existence of a
free disk containing two elements of $X$ would give rise to two
translation arcs not isotopic relative to $X$.  The second bulleted
item follows from the fact that $\mathbb R \times \{0\}$ and $\mathbb
R \times \{1\}$ are disjointly embedded lines.

As above, the proof of (3) will use the fact that $f_\#(\beta_i) \cap
\beta_i$ is exactly the terminal endpoint of $\beta_i$ but not the
fact that $\beta_i$ is isotopic to a translation arc.  It suffices to
assume that $f = Tr$, $x_1 = (0,0)$ and $x_2 = (\frac{1}{2}, 0)$.  We
may assume that $\rone \times \{0\}$ is a union of geodesics.  Let $K$
be any disk that contains $[0,1] \times \{0\}$ and that has geodesic
boundary.  We may assume that $Tr(\beta_1)$ is geodesic and in
particular that $Tr(\beta_1)$ and $\beta_1$ have no interior
intersections. In particular, if $b$ is a subpath of $\beta_1$ that intersects $\rone \times \{0\}$ exactly in its endpoints  then the
endpoints of $b$ span an interval on $\rone \times \{0\}$ with length
$\le 1.$  It follows that the union of all $[k,k+1]\times
\{0\}$'s that intersect $\beta_1$ is connected. 
Since $\beta_1 \cap [0,1] \times \{0\} = \emptyset$, $\beta_1$ is
disjoint from $[0,\infty)\times \{0\}$ and so can be isotoped into the
left half plane.  The analogous argument shows that $\beta_2$ can be
isotoped into the right half plane.  As geodesics have minimal
intersections in their isotopy classes, this completes the proof.\qed

\vspace{.1in}
We only use the following result for $k=2$ but the proof is the same for general $k$.

\begin{defn}  For all $k \ge 1$ let $R_k \co  \rtwo \to \rtwo$ be a homeomorphism that preserves $\rone \times \{0,1, \dots, k\}$ and that agrees with translation  to the left  by one unit on $\rone \times \{j\}$ for odd $1 \le j \le k$  and with translation to the right  by one unit on $\rone \times \{j\}$ for even $0 \le j \le k$.  If $f \co  \rtwo \to \rtwo$ is an orientation preserving homeomorphism and if $X = \cup_{i=0}^k\calo(x_i)$ for  $x_i\in \rtwo$ with distinct orbits then we say that  that {\it $f$ is a multiple Reeb class relative to $X$}  if the isotopy class of $f$ relative to $X$ is conjugate to the isotopy class of $R_k$ relative to $\cup_{j=0}^k \calo((0,j))$.
\end{defn} 

Disjoint subsets $Y$ and $Z$ of a circle $C$ {\it alternate around $C$} if any interval in $C$ with distinct endpoints in $Y$ [respectively $Z$] contains an element of $Z$ [respectively $Y$].

\begin{lemma} \label{multiple Reeb} Suppose that $f \co  D^2 \to D^2$ is an orientation preserving homeomorphism of the disk whose restriction to $\Int(D^2$) is fixed point free.  Suppose also that $\{x_0,\dots x_k\} \subset \Int(D^2)$, $k \ge 1$, satisfy:
\begin{description}
\item [(1)] the $\alpha$ and $\omega$ limit sets are single points $\alpha(x_i), \omega(x_i) \in \partial D^2$.  The $\alpha(x_i)$'s and $\omega(x_i)$'s are all distinct. 
\item [(2)] the elements of $\{\alpha(x_i)\}$ and $\{\omega(x_i)\}$
alternate around $\partial D^2$.
\end{description}
Then $f|_{\Int(D^2)}$ is a multiple Reeb class relative to $X = \cup_{i=0}^k\calo(x_i)$.
\end{lemma}

\noindent{\bf Proof of Lemma~\ref{multiple Reeb}}\qua   We will prove there are disjoint properly embedded lines $B_i \subset {\rm int}D^2$ that contain $\calo (x_i)$ and that are $f$-invariant up to isotopy relative to $X$.  The lemma then follows from (2).

By Theorem 2.7 of \cite{han:fpt} there is a line $\lambda \subset \Int(D^2)\setminus X$ such that:
\begin{itemize}
\item $\lambda$ is properly embedded in $\Int(D^2$).
\item $f(\lambda) \simeq \lambda$ relative to $X$.
\item both components $C_1$ and $C_2$ of $\Int(D^2) \setminus \lambda$ intersect $X$ non-trivially.
\end{itemize}

By (2), we may assume (see the last paragraph of the proof of Theorem
2.3 on page 258 of \cite{han:fpt}) that the ends of $\lambda$ converge
to single points in $\partial D^2 \setminus (\{\alpha(x_i)\} \cup
\{\omega(x_i)\}$.  After an isotopy relative to $X$ we may assume that
$f(\lambda) = \lambda$.

By induction on the number of orbits in $X \cap C_j$ , there are
disjointly embedded lines $B_i \subset {\rm int}(D^2)$ that contain
$\calo (x_i)$ and that are $f$-invariant up to isotopy relative to $X
\cap C_{j(i)}$ where $\calo (x_i) \subset C_{j(i)}$. There is no loss in
assuming that $B_i \subset C_{j(i)}$.  Since $B_i$ and $f(B_i)$ are
contained in $C_{j(i)}$, $B_i$ is invariant up to isotopy relative to
$X$ as desired. \qed

The following lemma is the analogue of Lemma~\ref{intersection}(2).  It's proof is virtually the same as that of Lemma~\ref{intersection}(2) and is left to the reader.

\begin{lemma} \label{mult Reeb recognition}  If $f$ is a multiple Reeb class relative to $X$ then no free disk contains two elements of $X$.
\end{lemma}

\section {Limit points on $\sinfty$}  
Suppose that $f\co M \to M$ is an orientation preserving homeomorphism of a hyperbolic surface and that $f$ is fixed point free and isotopic to the identity. We denote the identity lift by $\ti f \co  H \to H$.

 The following
proposition (which is essentially Proposition 5.1 of \cite{han:ze})
assigns two, not necessarily distinct, points in $\sinfty$ to each $\ti f$ orbit. 

\begin{prop} \label{single point} For any $\ti x \in H$ the $\alpha $ and $\omega$ limit sets  for $\ti f$ are single points $\alpha (\ti x)$ and $\omega(\ti x)$ in $\sinfty$.
\end{prop}

\begin{proof}
We give the argument for $\omega(\ti x)$. The Brouwer translation theorem implies that $\omega(\ti x) \subset \sinfty$.  Since $\ti f|_{\sinfty}$ is the identity, $\omega(\ti x)$ is connected.  If the lemma fails then  $\omega(\ti x)$ contains an interval $I \subset \sinfty$.  

Choose a geodesic $\ti \gamma \subset H$ that is the axis of a covering translation and that has endpoints in the interior of $I$;  let $U$ be the component of $H \setminus \ti \gamma$ whose closure intersects $\sinfty$ in a subinterval $J$ of $I$.  

We claim that there is a point $\ti y \in H$ whose backward orbit is
contained in the closure $\bar U$ of $U$.  To construct $\ti y$ choose
positive integers $m_i < n_i$ such that:
\begin{itemize}
\item $\ti f^{m_i}(\ti x), \ti f^{n_i}(\ti x) \not \in U$
\item $ \ti f^j(\ti x)  \in U$ for $m_i < j < n_i$
\item $n_i - m_i \to \infty$
\end{itemize}
Also choose a fundamental domain $W$ for the action of the covering
translation $T$ corresponding to $\ti \gamma$ and let $l_i$ satisfy
$\ti z_i =T^{l_i}(\ti f^{n_i-1}(\ti x)) \in W$.  Then $\ti z_i \in U \cap
W$ and $\ti f(\ti z_i) \not \in U$.  It follows that $\ti z_i$ lies
in a compact neighborhood of the intersection of $\ti \gamma$ with $W$
and hence that the accumulation set of the $\ti z_i$'s in $H$ is
non-empty.  Let $\ti y$ be any such accumulation point.  The second
and third item imply that $\calo ^-(\ti y) \subset \bar U$.

For the rest of the argument \tag s for $\ti f \co  H \to H$ are defined relative  to $X = \calo(\ti x) \cup \calo(\ti y)$.   

There are \tag s $\alpha_i$ for $\ti f^i(\ti x)$, $i>0$,   whose union is an immersed line with accumulation set equal to $I$. (For a proof that translation arcs can be locally constructed see, for example, Lemma 4.1  of \cite{han:fpt}). Similarly,  there are \tag s $\beta_j$ for $\ti f^{-j}(\ti y)$, $j > 0$,     whose union is an immersed line with accumulation set contained in $J$.  Since $J$ is contained in the interior of $I$ there are arbitrarily large $i$ and $j$ such that $\alpha_i \cap \beta_j \ne \emptyset$. This contradicts Theorem~\ref{fpt} and Lemma~\ref{intersection}.
\end{proof}

\begin{defn} \label{gamma}  If $\alpha(\ti x) \ne \omega(\ti x)$ then define $\ti \gamma(\ti x) \subset H$ to be the oriented geodesic connecting $\alpha(\ti x)$ to $\omega(\ti x)$ and define $\gamma(x) \subset S$ to be its projected image.  We say that {\it $\ti \gamma(\ti x)$ [respectively $\gamma(x)$] is the geodesic associated to $\calo(\ti x)$ [respectively $\calo(x)]$.} Denote the {\it union of all oriented geodesics in $M$ that occur as some $\gamma(x_i)$ by $\Gamma(f)$}.
\end{defn}

\section{Annuli}  

As in the previous section, $\ti f \co H \to H$ is the identity lift of a
fixed point free homeomorphism $f \co  M \to M$ of a connected surface with
negative Euler characteristic.
We state both of our main results, Lemma~\ref{no cycles} and
Lemma~\ref{near cycle}, before turning to their proofs.

\begin{defn} A {\it geodesic is birecurrent} if it is contained in the accumulation
set in the unit tangent bundle of its forward end and the accumulation set in the unit tangent bundle of its backward end.
\end{defn}

\begin{lemma} \label{no cycles} 
If $\Per(f)= \emptyset$ then:
\begin{description}
\item [(1)] Each $\gamma \in \Gamma(f)$ is simple.
\item [(2)] Birecurrent elements of  $\Gamma(f)$  have no transverse intersections. 
\end{description}
\end{lemma}

\begin{defn} 
Suppose that $U$ is a free disk for $f$ that contains both $ x$ and $
f^m( x)$ for some $x \in M$ and $m > 0$. Choose lifts $\ti x \in \ti U
\subset H$. The covering translation $T \co  H \to H$ satisfying $\ti
f^m(\ti x) \in T(\ti U)$ is called a {\it \nc} for $\ti x$ of {\it
period $m$}.  If the period is not relevant to the discussion then we
will simply say that $T$ is a near cycle for $\ti x$.
\end{defn}

\begin{remark} This definition of near cycle is less general than the one in \cite{han:ze} but is sufficient for our purposes.  
\end{remark}

\begin{lemma} \label{near cycle} If $\Per(f)= \emptyset$ and  $T \co  H \to H$  is a \nc\ for $\ti x$, 
 then $T^+$ and $T^-$ can not both lie in the same component of
$\sinfty \setminus(\{\alpha(\ti x),\omega(\ti x)\})$.
\end{lemma}

The proofs are technical so we begin with some motivation for our approach. Consider the situation of Lemma~\ref{no cycles}-(1).  We must show that if the geodesic $\gamma$ associated to some birecurrent point $x \in M$ has self intersections then $f$ has periodic points.  We have no control over the topology of  $\calo(x)$, it may even be dense, and so can not hope to analyze the isotopy class of $f$ relative to (the closure) of $\calo(x)$. Lifting to covers simplifies the relative isotopy class.  If we lift all the way to the universal cover $\ti f \co H \to H$, then  $\ti f$ is a translation class relative to $\calo(\ti x)$, so we completely understand the relative isotopy class,  but    the geodesic $\ti \gamma$ associated to $\ti x$ no longer has self intersections; we have oversimplified.

      The right cover to look is a cyclic one as in Lemma~\ref{canonical cyclic lift} with underlying space  an open annulus $\AO$.  If $\gamma$ has self intersections then there is a lift $\hat f \co  \AO \to \AO$ of $f \co  M \to M$ and a lift $\hat x$ of $x$ such that the geodesic $\hat \gamma$ associated to $\calo(\hat x)$ is properly immersed with a finite but non-zero number of self intersections.  Moreover the restriction of $\hat f$ to a neighborhood of the ends of $\calo(\hat x)$ is conjugate to the restriction of $\ti f$ to a neighborhood of the ends of $\calo(\ti x)$.  Thus, up to relative isotopy,  $\hat f$ acts on the ends of $\calo(\hat x)$ as a translation.  This is made precise in Lemma~\ref{half streamlines}.

    We are now in a situation just like the one we encountered in the proof of Theorem~\ref{canonical form}.  In both cases there is subsurface that contains all but finitely many of the punctures and on which we completely  understand the relative isotopy class.  In the case of Theorem~\ref{canonical form}, the infinitely punctured subsurface was fully invariant and we simply applied the Thurston classification theorem to the finitely punctured complement.  In our present case the infinitely punctured subsurface  is not fully invariant, so   Thurston's theorem does not apply to the finitely punctured complement.  The \lq fitted family\rq\ techniques of \cite{han:fpt} were developed for just this sort of application and it is these techniques that we employ in the  proofs of Lemmas~\ref{no cycles} and \ref {near cycle}.

  In \cite{han:fpt} the underlying space is the interior of a disk and here the underlying space is the interior of an annulus.  In addition to recalling definitions from  section 5 of \cite{han:fpt} we will use results that are proven there  but not explicitly stated.  These appear as Lemmas~\ref{quote0}, ~\ref{quote1} and ~\ref{quote2} in this paper.  Rather than repeat the proofs from \cite{han:fpt} we  give detailed references.

We now turn to the proof of Lemmas~\ref{no cycles} and \ref {near cycle}.

For $T \co H \to H$ a non-trivial covering translation with (not
necessarily distinct) endpoints $T^{\pm}$, let $\AO$ be the quotient
space of $H$  by the action of $T$ and let $\hat f \co  {\AO} \to {\AO}$
be the homeomorphism induced from $\ti f$.   We suppress the dependence
of $\AO$ and $\hat f$ on $T$ to simplify notation.  Denote the
quotient space of $\sinfty \setminus \{T^{+}, T^{-}\}$ under the
action of $T$ by $\partial_h \AO$.  If $T$ is hyperbolic, then
$\partial_h \AO$ compactifies $\AO$ to a closed annulus as in
Lemma~\ref{canonical cyclic lift}.  If $T$ is parabolic then
$\partial_h \AO$ compactifies one end of $\AO$.  In either case, $\hat f$ extends over $\partial_h \AO$ by the identity.  Assuming that neither
$\alpha(\ti x)$ nor $\omega(\ti x)$ is $T^+$ or $T^-$, denote their
images in $\partial_h { \AO}$ by $\alpha(\hat x)$ and $\omega(\hat x)$
respectively.

We  work relative to a finite union $\hat X =  \cup_{i=1}^r\calo(\hat x_i) \subset \AO$ of $\hat f$-orbits. Fix once and for all a complete hyperbolic structure on $\AO \setminus \hat X$.

The next lemma asserts that up to isotopy relative to $\hat X$, $\hat f$ acts on the ends of the orbits of $\hat X$ by translation. More precisely, for each $\hat x_i$ there is a  ray $B_i^+$ [respectively $B_i^-$ ]  that is forward [backward] $\hat f$ invariant up to isotopy relative to $\hat X$ and that contains the forward [backward] end of $\calo(\hat x_i)$.  As usual, the relative isotopy class of $\hat f$ is reflected in the  map $\hat f_\#$ that $\hat f$ induces on geodesics.

\begin{lemma}  \label{half streamlines}  Assume that $\Per(\hat f|_ {\AO}) = \emptyset$ and that    $\hat X =  \cup_{i=1}^r\calo(\hat x_i)$ where the $\alpha(\hat x_i)$'s are distinct  and the $\omega(\hat x_i)$'s are distinct. Define \tag s with respect to $\hat X$.  
Then there exist \tag s $\hat \beta_i^{\pm}$ for elements of $\calo(\hat x_i)$ such that:
\begin{description}
\item [(1)]  $\hat B^+_i = \cup_{n=0}^{\infty}  \hat f_\#^n(\hat \beta_i^+)$  is an embedded ray that converges to $\omega(\hat x_i)$. 
\item [(2)]  $\hat B^-_i = \cup_{n=0}^{\infty}  \hat f_\#^{-n}(\hat \beta_i^-)$  is an embedded ray that converges to $\alpha(\hat x_i)$.
\item [(3)] the $\hat B^{\pm}_i$'s are all disjoint.
\end{description}
\end{lemma}

\noindent{\bf Proof of Lemma~\ref{half streamlines}}\qua Choose lifts $\ti x_i \in H$ of $ \hat x_i$.  The argument for (1) takes place in a disk neighborhood $U_i \subset H \cup \sinfty$ of $\omega(\ti x_i)$ such that $T^k(U_i) \cap U_i = \emptyset$ for all $ k \ne 0$.  An isotopy relative to $\ti X = \cup_{k=-\infty}^{\infty} \cup_{j=0}^r \calo(T^k(\ti x_j))$  with support in $U_i$  projects to an isotopy relative to $\hat X$  with support in the image $\hat U_i \subset {\AO}$ of $U_i \cap H$.  

There are two cases.  In the first we assume that $T^k(\alpha(\ti
x_j)) \ne \omega(\ti x_i)$ for all $j$ and $k$.  In this case we may
assume that $U_i \cap \ti X \subset \calo(\ti x_i)$.  There are \tag s\ $\ti
\beta(n,\ti x_i)$ for $\ti f^n(\ti x_i)$ that converge to $\omega(\ti
x_i)$ as $n \to \infty$.  For all sufficiently large $n$, say $n \ge
N$, $\ti f(\ti \beta(n,\ti x_i))$ and $\ti \beta(n+1,\ti x_i)$ are
contained in $U_i$.  Lemma~\ref{intersection}-(1) implies that $\ti
f(\ti \beta(n,\ti x_i))$ and $\ti \beta(n+1,\ti x_i)$ are isotopic in
$H$ relative to $U_i \cap \ti X$ and hence isotopic in $U_i$ relative to $\ti X$.
Thus $\ti f_\#(\ti \beta(n,\ti x_i))= \ti \beta(n+1,\ti x_i)$.  It
follows that $\ti B^+_i = \cup_{n=0}^{\infty} \ti f_\#^n(\ti \beta(N,
\ti x_i) \subset U_i$ is an embedded ray that converges to $\omega(\ti
x)$.  Let $\hat \beta_i^+$ be the image of $\ti \beta(N,\ti x_i)$ and
$\hat B^+_i$ be the image of $\ti B^+_i$.

In the second case there are $j$ and $k$ such that
$T^k(\alpha(\ti x_j)) = \omega(\ti x_i)$.  We may assume that $U_i \cap
\ti X \subset \calo^+(\ti x_i) \cup \calo^-(T^k(\ti x_j))$.  As in the
previous case there are \tag s\ $\ti \beta(n,\ti x_i)$ for $\ti
f^n(\ti x_i)$ that converge to $\omega(\ti x_i)$ as $n \to \infty$.
If $\ti f(\ti \beta(n,\ti x_i))$ and $\beta(n+1,\ti x_i)$ are
contained in $U_i$ but are not isotopic relative to $U_i \cap \ti X$ then for
some $m < 0$, $\beta(n+1,\ti x_i) \cup \ti f_\#(\beta(n,\ti x_i))$
separates $\ti f^m(T^k(\ti x_j))$ from its $\ti f$-image.  In that
case every \tag\ for $\ti f^m(T^k(\ti x_j))$ intersects either $\ti
\beta(n+1,\ti x_i)$ or $\ti f_\#(\beta(n,\ti x_i))$.  Theorem~\ref{fpt}
and Lemma~\ref{intersection} imply that this can not happen for
sufficiently large $n$, say $n>N$.  We may therefore define $\hat
\beta_i^+$ and $\hat B^+_i$ as in the previous case.

This completes the proof of (1).  Condition (2) follows by replacing $f$ with $f^{-1}$.  It is clear from the construction that (3) will be satisfied if $N$ is sufficiently large.  \qed

\vspace{.1in}

\begin{lemma}\label{quote0}  Suppose that $\hat f \co  \AO \to \AO$, $\hat X$ and $\hat B_i^{\pm}$ are as in Lemma~\ref{half streamlines}.  Suppose that $\gamma$ is a simple geodesic in $\AO \setminus \hat X$ that is properly embedded in $\AO$ and  that intersects some $\hat B_i^+$ and some $\hat B_j^-$ infinitely often.  Then $\hat f_\#(\gamma)$ and $\gamma$ have non-trivial transverse intersection.
\end{lemma}

\proof  This follows easily from the fact that $\hat B_i^+$ and $\hat B_j^-$ translate in opposite directions.  A detailed proof is given in the last paragraph of the proof of Theorem 2.3 of \cite{han:fpt}.
\endproof

We next divide $\AO \setminus \hat X$ into a finitely many infinitely punctured subsurfaces on which we understand the dynamics and one finitely punctured subsurface on which we do not. 

There is one infinitely punctured subsurface $V_i^{\pm}$ for each $\hat B^{\pm}_i$.  Denote $\hat B^+_i \cap \hat X $ by $ \hat X_i^+ $. Then $V_i^+$ is uniquely determined by the properties that $\partial V_i^+$ is a single geodesic line properly embedded in $\AO$ , that $ V_i^+$ contains $\hat f_\#(\hat B^+_i)$  and that $V_i^+ \cap \hat X =\hat f(\hat X_i^+)$. 
 
It is useful to think in terms of the following topological model.
Identify $\AO$ with $\rtwo \setminus \{0\}$, $\hat B_i^+$ with 
$[1,\infty) \times \{0\}$ and $\hat X_i^+$ with $\N \times \{0\}$.  We may
assume that all the other $\hat B_j^{\pm}$'s are far away from $\hat
B_i^+$. Then $V^+_i$ can be identified with the set of points with
distance at most $\frac{1}{2}$ from $[2, \infty) \times \{0\}$.  Up to
isotopy relative to $\hat X$, $\hat f|_{V^+_i}$ agrees with
translation to the right by one unit followed by a vertical
compression toward the $x$-axis.  Thus up to isotopy relative to $\hat
X$, $\hat f(V_i^+) \subset V_i^+$ and $\cap_{k = 0}^{\infty} \hat
f^k(V_i^+) = \emptyset$.

The $V_i^-$'s are defined analogously using $\hat f^{-1}$ and $\hat B_i^-$ instead of $\hat f$ and $\hat B_i^+$. 

For each $i$ and $k$, the line $\hat f_\#^k(\partial V_i^{\pm})$ separates $\AO$.  Denote the complementary component that contains $\hat f^k(\hat X_i^{\pm})$ by $ \hat f_\#^k(V_i^{\pm})$.  Then $ \hat f_\#(V_i^+) \subset V_i^+$ and $\cap_{k = 0}^{\infty} \hat f_\#^k(V_i^+) = \emptyset$.  Similarly, $ \hat f^{-1}_\#(V_i^-) \subset V_i^-$ and $\cap_{k = 0}^{\infty} \hat f_\#^{-k}(V_i^+) = \emptyset$.

The subsurface $W = \AO
\setminus ( \hat X \cup (\bigcup_{i=1}^r \hat V_i^{\pm}))$ is finitely punctured.  We write $\partial W = \partial_+W \cup \partial_-W$ where $\partial_{\pm}W =  \cup_{i=1}^r \partial V^{\pm}_i$.   Then  $\hat f_\#(\partial_+W) \cap W = \emptyset$ and $\partial_-W  \cap \hat f_\#(W)  = \emptyset$.    We say that $W$ is the {\it Brouwer subsurface determined by the $B_i^{\pm}$'s}.

Denote  $\AO
\setminus ( \hat X \cup (\bigcup_{i=1}^r \hat f_{\#}(\hat V_i^{\pm})))$  by $\hat f_\#(W)$.    By part (2) of Lemma 3.5 of \cite{han:fpt} there is a homeomorphism $\hat f'$, isotopic to $\hat f$, such that $\hat f'(L) =  \hat f_\#(L)$ for each component $L$ of $\partial W$.  In particular $\hat f'(W) = \hat f_\#(W)$.

We  recall some more definitions from \cite{han:fpt}.

Let $\RH$ be the set of non-trivial relative homotopy classes determined by embedded arcs $(\tau, \partial \tau) \subset (W, \partial_+W)$.  Given any collection $ {\cal T} $ of elements of $\RH$ we will associate another collection $\hat f_\#({\cal T} ) \cap W$ of elements of $\RH$.  We abuse notation slightly and sometime write ${\cal T} =\{t_i\}$ where each $t_i \in \RH$; we do not assume that the $t_i$'s are distinct and it is essential that we allow multiplicity to occur.

Choose a homeomorphism $\hat f'$ as above.  For any arc $\tau \subset W$ with endpoints on $\partial_+W$, $\hat f'(\tau)$ is an arc in $\hat f_\#(W)$ with endpoints on $\hat f(\partial_+W)$; in particular, $\hat f'(\tau) \cap \partial_-W = \emptyset$ and $\partial \hat f'(\tau) \cap W = \emptyset$.  Let $\hat f_\#(\tau) \subset \hat f_\#(W)$ be the geodesic arc that is isotopic rel endpoints to $ \hat f'(\tau)$.  The components $\tau_1,\ldots,\tau_r$ of $\hat f_\#(\tau) \cap W$ are arcs in $W$ with endpoints in $\partial_+W$.  If $[\tau]$ denotes the element of $\RH$ determined by $\tau$ then we define $\hat f_\#([\tau] \cap W)$ to be $\{[\tau_1],\ldots,[\tau_r]\}$. It is shown in \cite{han:fpt} (see pages 249 - 250) that $\hat f_\#([\tau] \cap W)$ is well defined. 

   More generally if ${\cal T} $ is a finite collection of elements of $\RH$ then we define $\hat f_\#({\cal T} ) \cap W = \cup_{t \in {\cal T} }(\hat f_\#(t) \cap W)$.  Note that $\hat f_\#(.) \cap W$ can be iterated.  Inductively define $\hat f^n_\#(\tau) \cap W = \hat f _\#^{n-1}(\hat f_\#(\tau) \cap W) \cap W$.

We say that a finite collection  ${\cal T} =\{t_i\}$ of distinct elements of $\RH$ is {\it fitted} if the $t_i$'s are represented by simple disjoint arcs and if each element $s_j \in \hat f_\#(t_i \cap W)$ satisfies $\pm s_j \in {\cal T} $ where $\pm s_j$ means $s_j$ or $s_j$ with its orientation reversed.   We say that $t \in \RH$ {\it disappears under iteration} if $\hat f^n_\#(t) \cap W = \emptyset$ for some $n > 0$.  We say that a fitted family ${\cal T}  =\{t_i\}$ disappears under iteration if each $t_i$ does.

\begin{lemma} \label{quote1} Assume that $\Per(\hat f|_ {\AO}) = \emptyset$ and that    $\hat X =  \cup_{i=1}^r\calo(\hat x_i)$ where the $\alpha(\hat x_i)$'s are distinct  and the $\omega(\hat x_i)$'s are distinct.   If ${\cal T}  \subset \RH$ is a fitted family that does not disappear under iteration then there exists $t \in {\cal T} $ and $l > 0$ such that $\hat f^l_\#(t) \cap W = \{t,s_1,\ldots,s_m\}$ where each $s_i$ disappears under iteration.
\end{lemma}

\proof  The first half (up through the bottom of page 252) of the proof of Theorem~ 5.5 of \cite{han:fpt} applies to our current context without change to prove this lemma.
\endproof

A non-peripheral  line $\hat \lambda \subset \AO \setminus \hat X$ is  a {\it reducing line for $\hat f$ relative to $\hat X$} if it is properly embedded in $\AO$ and   is $\hat f$-invariant up to isotopy relative to $\hat X$.    
The isotopy class of $\hat f$ relative to $\hat X$ is  {\it reducible} if there a reducing line for $\hat f$ relative to $\hat X$. 

\begin{lemma} \label{quote2}  
Suppose that $W$ is a Brouwer subsurface for $\hat f \co  \AO \to \AO$
relative to $\hat X$ as above.  If there exists $t \in \RH$ and $l >
0$ such that $\hat f^l_\#(t) \cap W = \{t,s_1,\ldots,s_m\}$ where each
$s_i$ disappears under iteration then isotopy class of $\hat f^l$
relative to $\hat X$ is reducible and the reducing line for $\hat
f^l$ relative to $\hat X$ can be chosen to be geodesic in a complete
hyperbolic structure on $\AO \setminus X$.
\end{lemma}

\proof The last two paragraphs in the proof of Lemma 6.4 of
\cite{han:fpt} applies to our current context without change to prove
this lemma.  \endproof
 
We now put together these results from \cite{han:fpt} to prove our key technical lemma.    Disjoint subsets   of $\partial_h \AO$ are said to {\it alternate around $\partial_h \AO$} if their intersections with each component of $\partial_h \AO$  alternate around that component.

\begin{proposition} \label{reducible}  Assume that $\Per(\hat f|_ {\AO}) = \emptyset$ and that    $\hat X =  \cup_{i=1}^r\calo(\hat x_i)$ where either $r =  1$ or  the $2r$ points  $\{\alpha(\hat x_i),\omega(\hat x_i): 1 \le i \le r\}$ are   distinct.  Then there exists $l > 0$ such that the isotopy class of $\hat f^l$ relative to $\hat X$ is reducible.  Moreover,  if $r=1$ or if the elements of $\{\alpha(\hat x_i)\}$ and $\{\omega(\hat x_i)\}$ alternate around $\partial_h \AO$ then the reducing line $\hat \lambda$ for $\hat f^l$ relative to $\hat X$ can be chosen so that its ends converge to single points in $\partial_h \AO$.
\end{proposition}

\noindent {\bf Proof of Proposition~\ref{reducible}}\qua
Let $B_i^{\pm}$ and $\beta_i^{\pm}$ be as in Lemma~\ref{half streamlines}. Define   $W$ to be the Brouwer subsurface determined by the $B_i^{\pm}$'s.

Let ${\cal T} $ consist of one copy of each element of $\RH$ that occurs as a component of $\hat f^n_\#(\hat \beta_1^-) \cap W$ for some $n > 0$. Since the $\hat f^n_\#(\hat \beta_1^-)$'s are disjoint and 
simple and since   $W$ is
finitely punctured, ${\cal T} $ is finite.  It is fitted by construction.

If ${\cal T}$ disappears under iteration then    $\hat B_1 = \cup_{n=-\infty}^{\infty} \hat f_\#^n(\hat \beta_1^-) \subset \AO$ is a properly embedded $f_\#$-invariant line and a reducing line $\hat \lambda$ is obtained by pushing $\hat B_1$ off of itself into a non-contractible   complementary component.  The ends of $\hat \lambda$ converge to $\alpha(\hat x_1)$ 
and $\omega(\hat x_1)$ respectively.  In this case $l =1$.

Suppose now that ${\cal T}$ does not disappear under iteration.  Lemmas~\ref{quote1} and  ~\ref{quote2} imply the existence of $l > 0$ and a reducing line $\lambda$ for $\hat f^l$ relative to $\hat X$ that is realized as a geodesic in $\AO \setminus \hat X$.   Lemma~\ref{quote0} implies that each end of $\lambda$ intersects at most one $B_i^+$ or $B_j^+$.  There is therefore no obstruction to isotoping the end of $\lambda$ relative to the union of the $B_i^{\pm}$ so that it converges to a single point in $\partial_h \AO$.
\endproof

\medskip
\noindent{\bf Proof of Lemma~\ref{no cycles}}\qua 
If $\ti \gamma(\ti x)$
projects to a geodesic $\gamma \in \Gamma(f)$ that is not simple, then
there is a covering translation $T$ such that $T(\ti \gamma(\ti x))$
has transverse intersections with $\ti \gamma(\ti x)$.  In particular, the pair $\{\alpha(T^n(\ti x)), \omega(T^n(\ti x))\}$ links the pair
$\{\alpha(T^{n+1}(\ti x)), \omega(T^{n+1}(\ti x))\}$ for all $n$. Moreover, the endpoints $\alpha(\ti x)$ and $\omega(\ti x)$ of $\ti \gamma(\ti x)$ are distinct from the endpoints of the axis of $T$ so we may apply
Proposition~\ref{reducible} to $T$ and $\hat X = \calo(\hat x)$ where
$\hat x$ is the projected image of $\ti x$.  The resulting line $\hat
\lambda$ lifts to a line $\ti \lambda \subset H$ whose ends converge
to single points on $\sinfty$. Up to isotopy relative to the full pre-image $\ti X = \{\calo(T^n(\ti x_1)\} \subset
H$ of $\hat X$ we may assume that $\ti \lambda$ is $\ti f^l$-invariant for some $l > 0$.  But then $\ti \lambda$ defines a non-trivial partition of $\ti X$ in contradiction to the above linking property.   This proves (1).

To prove (2) we assume that $\gamma_1 =\gamma(x_1)$ and $\gamma_2 = \gamma(x_2)$ are birecurrent and have transverse intersections  and we argue to a contradiction.   Let $\ti \gamma_i =\ti \gamma(\ti x_i)$ be a lift of $\gamma_i$.  If the endpoints of $\ti \gamma_i$ differ by a covering translation, denote the image of the axis of this covering translation by   $\beta_i \subset M$.   Since $\gamma_1$ and $\gamma_2$ are birecurrent  and have transverse intersections, there are
finite subarcs $ \gamma_0( x_i)$ of $\gamma( x_i)$ (where we view $\gamma(x_i)$ as a possibly periodic, infinite geodesic) such that $
\gamma_0( x_1) \gamma_0( x_2)$ is an essential closed curve whose associated geodesic $\beta
\subset M$  is not a multiple of either $\beta_1$ or
$\beta_2$.  Lift $\gamma_0(x_i)$ to  $\ti \gamma_0(\ti x_i) \subset \ti \gamma(\ti x_i)$.       After replacing $\ti \gamma(\ti x_2)$  with a translate, we may assume  that the terminal endpoint of $\ti \gamma_0(\ti x_1)$ equals the initial endpoint of $\ti \gamma_0(\ti x_2)$. Let $T$ be the covering translation that carries the initial endpoint of $\ti \gamma_0(\ti x_1)$ to the terminal endpoint of $\ti \gamma_0(\ti x_2)$. 

 Denote $\alpha(T^n(\ti x_i))$ by $\alpha_i^n$ and $\omega(T^n(\ti x_i))$ by $\omega^n_i$.  Then  $\{\alpha^n_2, \omega^n_2\}$ links both  $\{\alpha^n_1, \omega^n_1\}$  and $\{\alpha^{n+1}_1, \omega^{n+1}_1\}$ for all $n$.  We claim that the $\alpha_i^n$'s and $\omega_i^n$'s are all disjoint.  We will focus on $\omega_i^0$, the other cases being analogous.  Since $T$ has no periodic points other than its two fixed points $\omega^0_1  \ne \omega_1^n$ for $n \ne 0$.  Since the axis of $T$ projects to $\beta$, our choice of $\beta$ guarantees that  $\omega_1^0 \ne \alpha_1 ^n$ for any $n$.  The $T^j(\ti \gamma(\ti x_2))$'s are disjoint because $\gamma(x_2)$ is simple.  By construction, $T^j(\ti \gamma(\ti x_2))$ intersects the component of $H \setminus \ti \gamma(\ti x_2)$ that is disjoint from $\omega^0_1$ and so is entirely contained in this component for $j \ne 0$.  In particular, $\omega_1^0$  is not equal to $\alpha_2^j$ or   $\omega_2^j$ for $j \ne 0$.  We already know that $\omega_1^0$  is not equal to $\alpha_2^0$ or   $\omega_2^0$ so we have verified the claim.

 Define $\AO$ using $T$ and let $\hat x_i$ be the image of $\ti x_i$.  By our verified claim  the $\alpha$ and $\omega$ limit sets of $\hat x_1$ and $\hat x_2$ are four distinct points in $\partial_h \AO$ and we may apply Proposition~\ref{reducible}.  The proof now concludes as in the previous case.   \qed

\medskip
\noindent{\bf Proof of Lemma~\ref{near cycle}}\qua Let $\ti x \in H$ be a
lift of $x \in M$.  Given a near cycle $T$ of $H$, if either
$\alpha(\ti x)$ or $\omega(\ti x)$ is $T^{\pm}$ there is nothing to
prove.  Assume then that  $\{\alpha(\ti x), \omega(\ti x)\} \cap \{T^{\pm}\} =
\emptyset$ and let $\hat f \co  {\AO} \to {\AO}$, $\hat X = \calo(\hat x)$, $l$
and $\hat \lambda$ be as in Proposition~\ref{reducible}.  As in the
proof of Lemma~\ref{no cycles}, there is a lift $\ti \lambda$ that
separates $\calo(\ti x)$ from $\calo(T(\ti x))$.  After an isotopy
relative to $\calo(\ti x) \cup \calo(T(\ti x))$ we may assume that
$\ti f^l$ preserves $\ti \lambda$.  

We claim that if $\{\alpha(\ti x),\omega(\ti
x)\}$ is contained in a single component of $\sinfty
\setminus\{T^{\pm}\}$ then $\ti f$ is Reeb class relative to $\calo(\ti
x) \cup \calo(T \ti x).$  By Theorem~\ref{fpt}, it suffices to show that $\ti f$ is not a translation class   relative to $\calo(\ti
x) \cup \calo(T \ti x)$ and therefore it suffices to show that $\ti f^l$ is not a translation class   relative to $\calo(\ti f^l,\ti
x) \cup \calo(\ti f^l,T \ti x).$   This follows from the observation that if $\{\alpha(\ti x),\omega(\ti
x)\})$ is contained in a single component of $\sinfty
\setminus\{T^{\pm}\}$ then $\ti \gamma(\ti x)$ and $T \ti \gamma(\ti
x)$ are anti-parallel.  

On the other hand, $T$ is a near cycle so
there is a free disk containing both $T(\ti x)$ and some $\ti f^t(\ti
x)$.  This contradicts Lemma~\ref{intersection}(2). \qed

\section{Birecurrent points}

As in the previous section, $\ti f \co H \to H$ is the identity lift of a fixed point free homeomorphism $f \co  M \to M$ of a hyperbolic surface.  We now also assume that  $\Per( f) = \emptyset.$  Let $\calb \subset M$ be
the set of birecurrent points of $f$ and let $\ti \calb \subset H$ be the
full pre-image of $\calb$.

The following definition is necessary to account for isolated punctures in $M$.

\begin{defn}  \label{def: rotates about}
We say that $\calo (x)$ {\it rotates about the isolated puncture $c$}
if for some (and hence all) lifts $\ti x \in H$ there is a parabolic
covering translation $T$ whose fixed point $P$ projects to $c$ such
that every near cycle for every point in $\calo(\ti x)$ is a positive
iterate of $T$.
\end{defn}

\begin{lemma}  \label{gamma exists} 
If $\ti x \in \ti \calb$ and $\ti \alpha(\ti x) = \ti \omega(\ti x) =
P$, then $P$ projects to an isolated puncture $c$ in $M$ and $\calo
(x)$ rotates about $c$.
\end{lemma}

\begin{proof}Suppose that $\ti
\alpha(\ti x) = \ti \omega(\ti x) = P$.  Lemma~\ref{near cycle}
implies that $P$ is $T^+$ or $T^-$ (or both) for an indivisible
covering translation $T$ and that all near cycles for elements of
$\calo (\ti x)$ are iterates of $T$.  

We recall from \cite{Fr2} that a {\it disk chain\/} for  $\ti f$ is an ordered set $(U_0,U_1,\dots,U_n)$ of embedded free
disks in the surface satisfying:

\begin{enumerate}
\item If $i\ne j$ then either $U_i=U_j$ or $U_i\cap U_j=\emptyset$

\item For $1\le i\le n$, there exists $m_i>0$ with $\ti f^{m_i}(U_i)\cap 
U_{i+1}\ne\emptyset.$
\end{enumerate}

If $U_0=U_n$ we say that $(U_0,\dots,U_n)$ is a {\it periodic disk chain}.
Proposition (1.3) of \cite{Fr2} asserts that if  an orientation
preserving homeomorphism of $\rtwo$  possesses a periodic disk
chain then it has a fixed point.  Thus $\ti f$ has no periodic disk chains.

For any disk chain $D=  (U_0,\dots,U_n)$, we denote the disk chain $(TU_0,\dots,TU_n)$ by $TD$.  If  $D_1 = (U_0,\dots,U_n)$ and $D_2 = (U_n,\dots,U_m)$ where all the $U_i$'s are disjoint or equal, then we denote the disk chain $(U_0,\dots,U_m)$ by $D_1 \cdot D_2$.

Suppose that both positive and negative iterates of $T$ occur as near
cycles for elements of $\calo (\ti x)$.   
Then there exist free disks $\ti U_1, \ti U_2$, elements $\ti x_i \in \ti U_i \cap \calo (\ti x)$ and  integers $p>0$ and $q < 0$ such that $D_1 = (\ti U_1, T^p(\ti U_1))$ and $D_2=(\ti U_2, T^q(\ti U_2))$ are disk chains.  Since all near cycles for elements of $\calo (\ti x)$ are iterates of $T$, every translate of $\ti U_i$ that intersect $\calo (\ti x)$ equals $T^a(\ti U_i)$ for some integer $a$.  If it is possible to choose $\ti U_1$ and $\ti U_2$ so that they have the same projection in $M$, do so.  If this is not possible, then birecurrence allows us to  replace $\ti U_i$ with an arbitrarily small disk neighborhood of $\ti x_i$.    We may therefore assume without loss that the projections $U_i \subset M$ of $\ti U_i$  are either disjoint or equal.   Birecurrence implies that there exist integers $r$ and $s$ such that the forward orbit of $\ti x_1$ intersects $T^r \ti U_2$ and the forward orbit of $\ti x_2$ intersects $T^s \ti U_1$.  In particular, $D_3 = (\ti U_1, T^r \ti U_2)$ and $D_4 = (\ti U_2, T^s \ti U_1)$ are disk chains.  For $i > -\frac{r+s}{q}$ the disk chain $D_5 = (D_3, \cdot T^r D_2 \cdot T^{r+q}D_2 \cdots T^{r+(i-1)q}D_2\cdot T^{r+iq}D_4)$ begins at $\ti U_1$ and ends at $T^a\ti U_1$ with $a < 0$.  But then $D_1$ and $D_5$ can be used to produce a disk chain that begins and ends at $\ti U_1$ in contradiction to the fact that there are no periodic disk chains for $\ti f$.     This contradiction proves that all near
cycles for $\calo (\ti x)$ are positive iterates of $T$ or all are
negative iterates of $T$.  After replacing $T$ by its inverse if
necessary, we may assume that all near cycles are positive iterates of
$T$.

If $T$ is parabolic, then we're done so suppose that $T$ is
hyperbolic.  Thus $T^+ = P$ and $T^- \ne P$.  Let us now switch our
focus to $\ti f^{-1}$.  All of its near cycles for $\ti x$ are
negative iterates of $T$.  But this contradicts the fact that $\ti x$
converges to $P \ne T^-$ under iteration by $\ti f^{-1}$.
\end{proof}

We will need the following facts about near cycles.

\begin{lemma} \label{delta} Suppose that $T$ is a near cycle for $\ti x$
of period $m$.   
\begin{description}
\item [(1)]  There is a neighborhood $\ti U$ of $\ti x$ such that $T$ is a
near cycle for $\ti y$ of period $m$ for each $\ti y \in \ti U$.
\item[(2)]   For any covering translation $S$, $STS^{-1}$ is a  near cycle
for $S\ti x$ of period $m$. 
\item [(3)]  There is a neighborhood $ \ti U$ of $\ti x$ so that if $\ti
f^{s}(x) \in S(\ti U)$ for some covering translation $S$,  then $STS^{-1}$
is a  near cycle for $\ti f^s(\ti x)$ of period $m$. 
\end{description}
\end{lemma}

\begin{proof} Items {1} and (2) follow from
the definition of near cycle.   If $\ti U$ is as in (1) and  $\ti
f^{s}(x) \in S(\ti U)$ then $T$ is a near cycle of period $m$ for $\ti
f^{s}(S^{-1}x)$.  Item (3) now follows from (2). 
\end{proof}

\begin{lemma} \label{small axes} If $\ti x \in \ti \calb$ and $\alpha(\ti
x) \ne \omega(\ti x)$ then one of the following holds.  
\begin{description}
\item [(1)]  Every near cycle for every point in $\calo (\ti x)$ has axis
$\ti \gamma(\ti x)$.
\item [(2)]  For every neighborhood $V$ of $\ti \alpha(\ti x)$ or $\ti
\omega(\ti x)$ there is a near cycle for an element of $\calo(\ti x)$
whose fixed points on $\sinfty$ are contained in $V$. 
\end{description}
\end{lemma}

\begin{proof}  We will prove  either that
(1) holds or that (2) holds for neighborhoods of $\omega(\ti x)$. Applying this
with $f^{-1}$ replacing $f$ completes the proof.  

 If $\omega(\ti x)$ is
the fixed point of a parabolic near cycle for an element of $\calo(\ti x)$
then (2) holds for neighborhoods of $\omega(\ti x)$ and we are
done.  Suppose then that  $\omega(\ti x)$ is not the fixed point of a
parabolic near cycle for an element of $\calo(\ti x)$ .  If there is a
parabolic near cycle for an element of $\calo(\ti x)$ , then
Lemma~\ref{near cycle} implies that its fixed point must be $\alpha(\ti
x)$.

For $t > 0$, define  $N_t =\{ k \in \mathbb Z :$ there is a near cycle for $ \ti f^k(\ti x)$ with period at most $ t\}$.  Fix $t > 0$ such that    $N_t \ne
\emptyset$.  Lemma~\ref{delta} and birecurrence implies that     $N_t$ is
bi-infinite.  For each $k \in N_T$ choose a near cycle  $T_k$ for  $ \ti f^k(\ti x)$with period at most $ t$.  Denote the axis of $T_k$ by $A_k$.  

   There is a uniform bound to the distance that  $T_k$ moves $\ti f^k(\ti x)$.  This has several consequences.   First, $\lim_{k \to
\infty} T_k(\omega(\ti x)) = \omega(\ti x)$.  This implies that $T_k$ is
hyperbolic for all large $k$ and hence by Lemma~\ref{delta}, hyperbolic for all $k$.
 Additionally, both the translation length of $T_k$ and the  distance between $\ti f^k(\ti x_k)$ and
$A_k$ are    uniformly bounded.  The set of axes of covering translations with uniformly bounded
translation length is discrete and closed in $H$.   Thus either there is
sequence $k_i \in N_t$ tending to infinity such that $A_{k_i}\to
\omega(\ti x)$ or there is a covering translation $T$ whose  axis $A$ has $\omega(\ti x)$ as an endpoint and satisfies
$A_k = A$ for all sufficiently large $k \in N_t$. If the former case holds
we are done, so assume that $A_k = A$ for all sufficiently large $k$, say $k > K$.  

 Choose any $k_1 < k_2  \in N_t$.  For any $\epsilon >
0$,   there exists $k > K$ and a covering translation $S$ such
that $\ti f^{k+j}(\ti x)$ is $\epsilon$ close to $\ti f^{k_1 + j}(S\ti x)$ for $0
\le j \le k_2-k_1 + t$.  In other words, the initial  $k_2-k_1+t$ points
in the orbits of $\ti f^k(\ti x)$ and of $\ti f^{k_1}(S\ti x)$ are within
$\epsilon$ of each other.  For sufficiently small $\epsilon$,
Lemma~\ref{delta} implies that $k,k+(k_2-k_1) \in N_t$, that $SA_k =
A_{k_1}$ and that $SA_{k+(k_2-k_1)} = S A_{k_2}$.  Since $A_k = A_{k+(k_2-k_1)}$,  $A_{k_1}= A_{k_2}$.  

We have now shown that $A_k =  A$ is independent of $k$.  The same argument that established  $\lim_{k \to
\infty} T_k(\omega(\ti x)) = \omega(\ti x)$ also proves that $\lim_{k \to
-\infty} T_k(\alpha(\ti x)) = \alpha(\ti x)$.  Thus $\alpha(\ti x)$ is an endpoint of $A$.    Since the choice of $T_k$ was arbitrary we have verified (1) for fixed $t$.  Since  $t$ was arbitrary, we have verified (1).
\end{proof}

\begin{defn}  \label{cross section}
A simple closed geodesic is a {\it partial cross-section} to
$\Gamma(f)$ if it is only crossed in one direction by elements of
$\Gamma(f)$ and if there is at least one such crossing.
 \end{defn}

\begin{lemma} \label{continuity}
 \begin{description}
\item [(1)] The map that assigns the unordered pair $\{\ti \alpha(\ti x), \ti \omega(\ti x)\}$ to $\ti x \in \ti \calb$ is continuous.
\item [(2)]  If $x \in \calb$ and $\gamma(x)$ is defined then $\gamma(x)$ is birecurrent.
\item [(3)]  If $x \in \calb$ and  $\gamma(x)$ is defined but is not a simple closed curve,  then $\gamma(x)$ intersects a partial cross-section $\alpha$. In this case, the set $\{y \in \calb: \gamma(y)$ crosses $\alpha\}$ is open in $\calb$. \end{description}
\end{lemma}

\begin{proof}
We first show that if $\gamma(x)$ is defined, then continuity at $\ti x$ implies that $\gamma(x)$ satisfies (2) and (3).

Choose $ x_i =  f^{n_i}( x)$ with $n_i > 0$ such that $ x_i \to x$; choose also lifts $\ti x_i \to \ti x$. The geodesics $\ti \gamma(\ti x_i)$ are translates of $\ti \gamma(\ti x)$.  If $\ti \gamma (\ti x)$ is the axis of a covering translation, then $\gamma(x)$ is a simple closed curve and we are done. We may therefore assume that the $\ti \gamma(\ti x_i))$'s are distinct as unoriented lines.  We may also assume
there is a free disk that contains $ x$ and each $x_i$. 

The continuity assumption implies that $\ti \gamma (\ti x_i) \to \ti
\gamma(\ti x)$ as unoriented lines.  From Lemma~\ref{no cycles}(1) it
follows that the $\ti \gamma (\ti x_i)$'s are disjoint.  Lemmas
~\ref{multiple Reeb} and \ref{mult Reeb recognition} therefore imply
that for sufficiently large $i$, say $i > N$, the $\ti \gamma(\ti
x_i)$'s are consistently oriented (i.e.\ the lines are parallel and not
anti-parallel).  Let $T_i$ be the covering translation that carries
$\ti f^{n_i}(\ti x)$ to $\ti x_i$ and hence also $\ti \gamma(\ti x)$
to $\ti \gamma(\ti x_i)$.  Choose $\ti y \in \gamma(\ti x)$.  Since $\ti f^{n_i}(\ti x) \to \omega(\ti
x)$ and $\ti \gamma(\ti x_i)\to \ti \gamma(\ti x)$ there are points
$\ti y_i \in \ti \gamma(\ti x)$ such that $\ti y_i \to \omega(\ti x)$
and $T_i(\ti y_i) \to \ti y$.  This proves that the accumulation set of the  forward end of $\gamma(x)$ contains $\ti y$ and hence all of $\gamma(x)$ since $\ti y$ was arbitrary.    By symmetry, $\gamma(x)$ 
is also contained in the accumulation set of its backward end so
$\gamma(x)$ is birecurrent and we have verified (2).

    Lemma~\ref{no cycles}(2) now implies that $\gamma(x)$ does not intersect any other $\gamma(y)$.  Lemmas
~\ref{multiple Reeb} and \ref{mult Reeb recognition} therefore imply that if $\ti \gamma(\ti y)$ is contained in the region bounded by $\ti \gamma(\ti
x_i)$ and $\ti \gamma(\ti x_j)$ for $i,j > N$, then the orientation on
$\ti \gamma(\ti y)$ agrees with the orientation $\ti \gamma(\ti x_i)$
and $\ti \gamma(\ti x_j)$.  Let $\ti \alpha_0$ be a geodesic arc connecting $\ti \gamma(\ti
x_i)$ and $\ti \gamma(\ti x_j)$ and let $\alpha_0$ be its projected image.  Then elements of $\Gamma(f)$ cross $\alpha_0$ only in one direction.  Extend $\alpha_0$ to a simple closed curve $\alpha'$ by following $\gamma(x)$ from one intersection with $\alpha_0$ to the next.  The geodesic $\alpha$ corresponding to $\alpha'$ is a partial cross-section.  This establishes the first part of (3).  The second part follows from the fact that a pair of geodesics in $H$ have non-empty transverse intersection if and only if their endpoints link in $\sinfty$ and the fact that linking a particular pair of points in $\sinfty$ is an open condition on pairs of points  in $\sinfty$.

It remains to establish (1).  Suppose that $\ti x_i, \ti x \in \ti
\calb$ and that $\ti x_i \to \ti x$.  We assume at first that $\alpha(\ti x) \ne \omega(\ti x)$.   According to Lemma~\ref{small axes}
there are two cases to consider.  In the
first, every neighborhood $U^+$ of $\omega(\ti x)$ contains $T_1^+$ and $T_1^-$ where $T_1$ is
 a near cycle for an element of $\calo (\ti x)$ and every
neighborhood $U^-$ of $\alpha(\ti x)$ contains $T_2^+$ and $T_2^-$ where $T_2$ is
 a near cycle  for an element of $\calo (\ti x)$.  For sufficiently large
$i$ both $T_1$ and $T_2$ are near cycles for elements of $\calo (\ti
x_i)$. Corollary~\ref{near cycle} implies that $\{\alpha(\ti x_i),
\omega(\ti x_i)\}$ intersects both $U^+$ and $U^-$ for all
sufficiently large $i$. This completes the proof of (1) in the first
case and so also the proof of (2) and (3)  in the first case.

In the second case,  Lemma~\ref{small axes} implies that  $\ti \gamma(\ti x)$ is the axis of a near cycle for $\ti x$ so we have proved (2) and (3).  As in the first case, for all sufficiently large $i$, say $i > N$, $\ti \gamma(\ti x)$ is the axis of a near cycle for an element of $\calo (\ti x_i)$.  Lemma~\ref{no cycles} and (2) imply that $\ti \gamma(\ti x_i)$ and $\ti \gamma(\ti x)$ are either disjoint or equal for  $i > N$. Lemma~\ref{near cycle} therefore implies that they are equal. 

    It remains to establish (1) when $\ti \alpha(\ti x) = \ti \omega(\ti x)$.  By Lemma~\ref{gamma exists},  $\ti \alpha(\ti x) = \ti \omega(\ti x) = T^{\pm}$ for a parabolic covering translation and  every near cycle for $\calo (\ti x)$ is an iterate of $T$. For sufficiently large
$i$, iterates of $T$ are near cycles for elements of $\calo (\ti
x_i)$. By Lemma~\ref{near cycle} either $\alpha(\ti x_i)$ or $\omega(\ti x_i)$ is $T^{\pm}$.  This rules out the possibility that $\gamma(x_i)$ is defined and birecurrent.  The only other possibility is that $\ti \alpha(\ti x_i) = \ti \omega(\ti x_i) = T^{\pm}$.
\end{proof}

\section{Proof of Theorem~\ref{thm: periodic point}}

The following result is implicit in the paper of Atkinson \cite{A}.
Since it is not explicit there we give the (short) proof here.  This
proof is essentially the same as an argument in \cite{A}.

\begin{prop}\label{prop: atkinson}
Suppose $T\co  X \to X$ is an ergodic automorphism of a probability space
$(X,\nu)$ and let $\phi\co  X \to \R$ be an integrable function with
$\int \phi \ d\nu = 0.$ Let $S(n,x) = \sum_{i=0}^{n-1} \phi( T^i(x))$.
Then for any $\varepsilon >0$ the set of $x$ such that $|S(n,x)| <
\varepsilon$ for infinitely many $n$ is a full measure subset of $X$.
\end{prop}

\begin{proof}
Let $A$ denote the set of $x$ such that $|S(n,x)| < \varepsilon$ for
only finitely many $n$.  We will show the assumption $\mu(A) > 0$
leads to a contradiction.  Suppose $\mu(A) > 0$ and let $A_m$ denote
the subset of $A$ such that $|S(i,x)| < \varepsilon$ for $m$ or fewer
values of $i$.  Then $A = \cup A_m$ and there is an $N >0$ such that
$\mu(A_N) > p$ for some $p >0.$

The ergodic theorem applied to the characteristic function of $A_N$
implies that for almost all $x$ and all sufficiently large $n$ (depending
on $x$) we have 
\[
\frac{card( A_N \cap \{T^i(x)\ |\ 0 \le i < n\})}{n} > p.
\]
We now fix an $x \in A_N$ with this property.
Let $B_n = \{i\ |\ 0 \le i \le n 
\text{ and } T^i(x) \in A_N\}$ and $r = card(B_n)$; then $r > np$.
Any interval in $\R$ of length $\varepsilon$ which 
contains $S(i,x)$ for some $i \in B_n$  contains at most $N$ values of
$\{S(j,x) : j > i\}.$
Hence any interval of length $\varepsilon$ contains at most 
$N$ elements of $\{ S(i,x)\ |\ i\in B_n\}.$
Consequently an interval containing 
the $r$ numbers $\{ S(i,x)\ |\ i \in B_n\}$ must have length at least
$r\varepsilon/N$.  Since $r > np$ this length is $> np\varepsilon/N.$
Therefore
\[
\sup_{0 \le i \le n} |S(i,x)| > \frac{np\varepsilon}{2N},
\]
and hence by the ergodic theorem, for almost all $x \in A_N$ 
\[
\Big | \int \phi\ d\mu \Big |
= \lim_{n \to \infty} \frac{|S(n,x)|}{n} 
= \limsup_{n \to \infty} \frac{|S(n,x)|}{n}
> \frac{p\varepsilon}{2N} > 0.
\]
This contradicts the hypothesis so our result is proven.
\end{proof}

\begin{cor}\label{cor: atkinson}
Suppose $T\co  X \to X$ is an automorphism of a Borel probability space
$(X,\mu)$ and $\phi\co  X \to \R$ is an integrable function.
Let $S(n,x) = \sum_{i=0}^{n-1} \phi( T^i(x))$ and suppose 
$\mu(P) > 0$ where $P = \{x \ |\ \lim_{n \to \infty} S(n,x) = \infty\}.$
Let 
\[
\hat \phi(x) = \lim_{n \to \infty} \frac{S(n,x)}{n}.
\]
Then $\int_P \hat \phi \ d\mu > 0.$  In particular $\hat \phi(x) >0$ for a set
of positive $\mu$-measure.
\end{cor}
\begin{proof}
By the ergodic decomposition theorem there is a measure $m$ on the
space $\M$ of all $T$ invariant ergodic Borel measures on $X$ 
with the property that for any $\mu$ integrable function 
$\psi \co  X \to \R$ we have 
$\int \psi \ d\mu = \int_\M I(\psi,\nu) \ dm$ where $\nu \in \M$
and $I(\psi,\nu) = \int \psi\ d\nu.$

The set $P$ is $T$ invariant. Replacing $\phi(x)$ with 
$\phi(x) \X_P(x),$ where $\X_P$ is the characteristic function
of $P,$ we may assume that $\phi$ vanishes outside $P$.
Then clearly $\hat \phi(x) \ge 0$ for all $x$ for which it exists.
Let $\M_P$ denote $\{ \nu \in \M\ |\ \nu(P) > 0 \}$.  
If $\nu \in \M_P$ the fact that $\hat \phi(x) \ge 0$ and the ergodic theorem
imply that $I(\phi,\nu) = \int \phi\ d\nu  = \int \hat \phi\ d\nu  \ge 0$. 
Also Proposition~\ref{prop: atkinson}
implies that $\int \phi\ d\nu  = 0$ is impossible so
$I(\phi,\nu) > 0.$
Then $\mu(P) =  \int I( \X_P, \nu)\ dm = \int \nu(P)\ dm = 
\int_{\M_P} \nu(P)\ dm.$   This implies $m(\M_P) > 0$ since
$\mu(P) > 0.$

Hence
\[
\int \hat \phi \ d\mu = \int \phi \ d\mu = \int I( \phi, \nu)\ dm \ge \int_{\M_P} I( \phi, \nu)\ dm
>0
\]
since $I( \phi, \nu) >0$ for $\nu \in \M_P$ and $m(\M_P) > 0.$
\end{proof}

We are now prepared to complete the proof of the following result from
\S1.

\medskip
\noindent
{\bf Theorem~\ref{thm: periodic point}}\qua
{\sl 
Suppose $F\co  S \to S$ is a non-trivial, Hamiltonian surface
diffeomorphism and that if $S = S^2$ then $\Fix(F)$ contains at least three
points. Then $F$ has periodic points of arbitrarily high period.
}
\medskip
\begin{proof}  We first observe that $F$ has infinite order.  If $R(F) \ne \emptyset$ then this follows from the fact that the isotopy class of $F$ relative to $\Fix(F)$ has infinite order. If $R(F)  = \emptyset$ this follows from Lemma~\ref{s2}.

Suppose that the periods of
points in $\Per(F)$ are bounded.  We will show that this
leads to a contradiction. 
 Replacing $F$ with an iterate
we may assume $\Fix(F) = \Per(F) \ne \emptyset$.  By \cite{brn-kis} and Proposition~\ref{prop: per point}, $M =S \setminus \Fix(F)$ is $F$ invariant and $f = F|_M \co  M \to M$ is isotopic to the identity.

 No component of $M$ can be a disk by Theorem~\ref{thm: Brouwer}.  We
first handle the case that some component $U$ is an annulus.  In this
case by Lemma~\ref{lem: extension} $U$ can be compactified to a closed
annulus $\bar U$ in such a way that $f$ extends to the identity on at
least one boundary component (this is where the assumption that
$\Fix(f)$ has at least three points if $S = S^2$ is used).  If the
lift of $f|_{\bar U}$ which fixes points on the boundary were to have
mean rotation number zero, $f|_{\bar U}$ would have an interior fixed
point by Theorem 2.1 of \cite{Fr6}, which is impossible.  Hence this
lift has non-zero mean rotation number and consequently a point of
non-zero rotation number.  Also $f|_{\bar U}$ is chain recurrent since
it preserves a finite measure positive on open sets.  It now follows
from Theorem~\ref{thm: Birkhoff} that $f$ has points of arbitrarily
high period.

We are left with the case that all components of $M$ have a 
hyperbolic structure.    If $M$ is not connected we replace it by one 
component.   Let $\ti f \co  H \cup \sinfty \to H \cup \sinfty$
be the identity lift of $f$.

We continue with the notation of the previous sections.  Recall that
by the Poincar\'e recurrence theorem the set of birecurrent points of
$M$ has full measure.  Since there are only countably many simple closed geodesics and isolated punctures in $M$, Lemma~\ref{gamma exists} and Lemma~\ref{continuity} imply that there is a set $P \subset S$ of positive measure such that one of the following conditions is satisfied by all $x \in P$: 
\begin{enumerate}
\item   $\gamma(x) = \gamma$ for some simple closed geodesic $\gamma$
\item  $\calo (x)$  rotates about $c$ for some isolated puncture $c$
\item  $\gamma(x)$ crosses $\alpha$ for some partial cross section $\alpha$. 
\end{enumerate}

We consider these possibilities in order.  

First assume $\gamma(x) = \gamma$ is a simple closed geodesic.  Let
$\ti \gamma$ be a lift of $\gamma$ and let $\ti P(\ti \gamma) =\{ \ti
x \in H: \ti x$ projects into $P$ and $\ti \gamma(\ti x) = \ti
\gamma\}$.  Recall that the canonical cyclic cover $\hat f \co  \A \to
\A$ associated to $\gamma$ is induced from $\ti f$ by taking the
quotient space of $H \cup \sinfty$ under the action of the indivisible
covering translation $T$ with axis equal to $\ti \gamma$.  By
Lemma~\ref{lem: per points f hat} it suffices to show that $\ti f$
considered as a lift of $\hat f$ has points of two different rotation
numbers.  Denote the image of $\ti P(\ti \gamma)$ in $\A$ by $\hat P$.
If $T_0$ is another covering translation, not a power of $T$, then
$\ti \gamma(T_0(\ti x)) = T_0(\ti \gamma) \ne \ti \gamma$ so the
covering map $\pi \co  \Int(\A) \to M$ restricts to a bijection between
$\hat P $ and $P$.  In particular there is a finite measure $\mu$ on
$\hat P$ defined by $\mu(A) = \omega(\pi(A))$.  The identity lift of
$f$ is a lift of $\hat f$.  We can now apply Corollary~\ref{cor: atkinson}
to $\hat f\co  \A \to \A$ with measure $\mu$ and function $\phi
= \delta,$ a homological displacement function, so that $\hat \phi(x)
= \rho(x, \ti f).$ It follows that there is a point of $\A$ with
positive rotation number.  Since any point of $\partial \A$ has zero
rotation number for $\hat f$, this completes the proof in the first
case.

Next we consider the case that
$\calo (x)$  rotates about $c$ for some isolated puncture $c$.
Recall from Definition~\ref{def: rotates about} that this means
if $x \in P$ then for a lift $\ti x \in H$ there is a parabolic
covering translation $T$ whose fixed point $p$ corresponds to $c$ such
that every near cycle for every point in $\calo(\ti x)$ is a positive
iterate of $T$.

Let $\ti P(p) =\{ \ti x \in H: \ti x$ projects into $P$ and $\ti
\omega(\ti x) = \ti \alpha(\ti x) = p\}$.  Let $\hat f \co  \A \to \A$ be
the canonical cyclic cover associated to $\gamma$ and recall that
$\hat f|_{{\rm int}\A}$ is induced from $\ti f|_H$ by taking the
quotient space under the action of $T$.  As in the previous case the
image $\hat P \subset \A$ of $\ti P(p)$ has a finite invariant measure
defined by $\mu(A) = \omega(\pi(A))$.  By Lemma~\ref{canonical cyclic
lift} and Lemma~\ref{lem: per points f hat} it suffices to show that
$\ti f$ considered as a lift of $\hat f$ has points of positive
rotation number.

 Let $\hat U \subset {\rm int} \A$ be a lift of an
open free disk $U \subset M$ for $f$ such that $\mu(\hat U_P) > 0$
where $\hat U_P = \hat U \cap \hat P.$ A
full measure set of $U_P = U \cap P$  is recurrent.  But the fact that $\calo
(x)$ rotates about $c$ for $x \in P$ means that if $x \in U_P \cap
f^k(U_P),$ for $k >0$, and $\hat x \in \hat U_P \cap \pi^{-1}(x)$ then
$\hat x \in \hat U_P \cap \hat f^k(\hat U_P).$  It follows that if
$\hat V = \cup_{n \ge 0} \hat f^n(\hat U_P)$ and 
$V = \cup_{n \ge 0} f^n( U_P)$ then $\pi|_{\hat V}\co  \hat V \to V$
is a measure preserving bijection on a set of full measure.  In
particular $\mu(\hat V)$ is finite.

Let $\delta\co  \A \to \R$ be a homological displacement function for
the lift $\hat f$ and let $\delta_U\co  \hat U_P \to \R$ be defined by
\[
\delta_U( x) = \sum_{m = 0}^{n-1} \delta( \hat f^m(x)).
\]
where $x \in \hat U_P$,\ $\hat f^n(x) \in \hat U_P$, and $\hat f^k(x) \notin
\hat U_P$ for $1 \le k < n.$ Note in this case that if $\ti U_P$ is a lift
of $\hat U_P$ and $\ti x \in \ti U_P$ is a lift of $x,$ then $\ti
f^n(\ti x) \in T^k(\ti U_P)$ for some $k >0.$ Hence if we identify
$H_1( \A, R)$ with $\R$ we can consider $\delta_U$ as a real valued
function which takes on only strictly positive integer values.  In
effect $\delta_U$ is a homological displacement function for the first
return map of $\hat f$ on $U_P$.

We claim that $\int_{\hat V} \delta\ d\mu = \int_{\hat U_p} \delta_U\ d\mu.$
This can be seen by setting 
\[
W_n = \{ x \in \hat U_P | \hat f^n(x) \in \hat U_P
\text{ and } \hat f^k(x) \notin \hat U_P \text{ for } 0 < k < n\}
\]
and letting $V_n = \cup_{m=0}^{n-1} \hat f^m(W_n).$ Then 
\[
\int_{W_n} \delta_U\ d\mu 
= \int_{W_n} \sum_{m = 0}^{n-1} \delta( \hat f^m(x))\ d\mu 
=  \sum_{m = 0}^{n-1} \int_{\hat f^m(W_n)}\delta\ d\mu 
= \int_{V_n}\delta\ d\mu.
\]
Since the $W_n$'s are pairwise disjoint as are the $V_n$'s and 
since $\hat U_P = \cup W_n$ and $\hat V = \cup V_n$ the claim 
$\int_{\hat V} \delta\ d\mu = \int_{\hat U_p} \delta_U\ d\mu$
follows.

Hence
\[
\int_{\hat V} \rho( x, \ti f)\ d\mu = \int_{\hat V} \delta\ d\mu
= \int_{\hat U_p} \delta_U\ d\mu >0,
\]
and it follows that $\rho( x, \ti f) > 0$ on a set of positive 
measure.  Since $\hat f$ is covered by the identity lift
it has fixed points in the boundary of $\A$ and hence it
has points of two distinct rotation 
numbers. By Lemma~\ref{lem: per points f hat} $\hat f$
has points of arbitrarily high period.

We are left with the case that $\gamma(x)$ crosses the partial cross section $\alpha$, 
which we will show contradicts the hypothesis
that the flux homomorphism is zero. We assume that all crossings of $\alpha$ by elements of $\Gamma(f)$ are in the positive direction.  Let $F_t$ be an isotopy of 
$S$ with $F_0 = id$ and $F_1 = F.$
Let $\delta\co  S \to H_1(S,R)$ be a homological displacement function.
Then
\[
\rho_\mu(F_t) \wedge [\alpha] = \int \rho(x,F_t) \wedge [\alpha]\ d\mu
= \int \delta(x) \wedge [\alpha]\ d\mu.
\]
Also since
there is no $x$ for which $\gamma(x)$ crosses $\alpha$ in the negative
sense we have $\rho(x, F_t)\wedge [\alpha] \ge 0$ for all $x$ in a
full measure subset of $S$ and hence
\[
\rho_\mu(F_t) \wedge [\alpha] = \int \rho(x,F_t) \wedge [\alpha]\ d\mu
 \ge  \int_P \rho(x,F_t) \wedge [\alpha]\ d\mu.
\]
We want now to apply Corollary~\ref{cor: atkinson} with $\phi(x)
= \delta(x)\wedge [\alpha],$ which will mean that
$\hat \phi(x) = \rho(x, F_t)\wedge [\alpha].$
If $S(n,x) = \sum_{i=0}^{n-1} \delta( F^i(x))\wedge [\alpha]$ the fact that
$\gamma(x) = \gamma$ for all $x \in P$ implies that $\displaystyle{
\lim_{n \to \infty} S(n, x)) = \infty}$ for all $x \in P$.  
Thus by Corollary~\ref{cor: atkinson} we conclude 
\[
\rho_\mu(F_t) \wedge [\alpha] \ge 
\int_P \rho(x, F_t)\wedge [\alpha] \ d\mu = \int_P \hat \phi\ d\mu >0.
\]
Therefore by Proposition~\ref{prop: mean = flux} $\F(F_t)([\alpha]) \ne
0$ contradicting the hypothesis that $F$ is Hamiltonian.
\end{proof}

\noindent{\bf Proof of Theorem~\ref{thm: non-isotopic}}\qua Replacing $F$ by a
power we may assume that $F$ has more than $2g+2$ fixed points
where $g$ is the genus of $S$.
Let $Q$ be a finite set consisting of more than $2g+2$ fixed points.
By Theorem~\ref{thm: periodic point}, we can choose a periodic point $x$ whose 
period is $p > 1.$  Letting $P = Q \cup orb(x)$ we observe that the Thurston
canonical form of $F$ relative to $P$ cannot be periodic.  This is because
any non-trivial periodic diffeomorphism of $S$ can have at most $2g+2$ fixed
points (see \cite{FF}).  

Hence some iterate of the Thurston
canonical form of $F$ relative to $P$  must be reducible with either a pseudo-Anosov component or
a non-trivial Dehn twist along a reducing curve.  In either of these
cases $F^p$ is not isotopic to the identity relative to $P$.
\qed

\medskip
\noindent{\bf Proof of Theorem~\ref{thm: all high periods}}\qua This is an immediate consequence of Theorem~\ref{thm: non-isotopic} and Proposition~\ref{prop: per point}. 

\rk{Acknowledgements}

John Franks is supported in part by NSF grant 
DMS0099640.\nl Michael Handel is supported in part by NSF grant DMS0103435.

\end{document}